\documentclass{amsart} 
\usepackage{amssymb,amsmath,latexsym,times,xcolor,hyperref,tikz,bm}

\hypersetup{colorlinks=true, linkcolor=blue, citecolor=blue}

\numberwithin{equation}{section}

\theoremstyle{plain}
\newtheorem{thm}{Theorem}[section]

\newtheorem{cor}[thm]{Corollary}

\newtheorem{lemma}[thm]{Lemma}

\theoremstyle{definition}
\newtheorem{example}{Example}

\DeclareMathOperator{\cl}{cl}
\DeclareMathOperator{\ch}{ch}
\DeclareMathOperator{\fl}{flag}
\DeclareMathOperator{\comp}{comp}

\title{Matroids with different configurations and the same
  $\mathcal{G}$-invariant}

\author{Joseph E.~Bonin} \address{Department of Mathematics\\ The
  George Washington University\\  Washington, D.C. 20052, USA}
\email{jbonin@gwu.edu}

\date{\today.}

\begin{document}

\begin{abstract}
  From the configuration of a matroid (which records the size and rank
  of the cyclic flats and the containments among them, but not the
  sets), one can compute several much-studied matroid invariants,
  including the Tutte polynomial and a newer, stronger invariant, the
  $\mathcal{G}$-invariant.  To gauge how much additional information
  the configuration contains compared to these invariants, it is of
  interest to have methods for constructing matroids with different
  configurations but the same $\mathcal{G}$-invariant.  We offer
  several such constructions along with tools for developing more.
\end{abstract}
    
\maketitle

\section{Introduction}

The configuration of a matroid, which Eberhardt \cite{config}
introduced, is obtained from its lattice of cyclic flats (that is,
flats that are unions of circuits) by recording the abstract lattice
structure along with just the size and rank of each cyclic flat, not
the set.  Eberhardt proved that from the configuration of a matroid
$M$, one can compute its \emph{Tutte polynomial}, 
  $$T(M;x,y)=\sum_{A\subseteq E(M)}(x-1)^{r(M)-r(A)}(y-1)^{|A|-r(A)}.$$
The data recorded in the Tutte polynomial is the multiset of size-rank
pairs, $(|A|,r(A))$, over all $A\subseteq E(M)$.  The Tutte polynomial
is one of the most extensively studied invariants of a matroid (see,
e.g., \cite{Tutte,handbook}).

Strengthening Eberhardt's result, Bonin and Kung \cite[Theorem
7.3]{catdata} showed that from the configuration of a matroid $M$, one
can compute its $\mathcal{G}$-invariant, $\mathcal{G}(M)$.  Derksen
\cite{G-inv} introduced the $\mathcal{G}$-invariant and showed that
the Tutte polynomial can be computed from it.  The perspective on the
$\mathcal{G}$-invariant that we use is the reformulation introduced in
\cite{catdata}: $\mathcal{G}(M)$ records the multiset of sequences of
sizes of the sets in flags (maximal chains of flats) of $M$. (Section
\ref{sec:back} gives a more precise formulation.)  This information
about flags just begins to suggest the wealth of information that
$\mathcal{G}(M)$ captures beyond what the Tutte polynomial contains;
other examples (from \cite[Section 5]{catdata}) include the number of
saturated chains of flats with specified sizes and ranks, the number
of circuits and cocircuits of each size (in particular, the number of
spanning circuits), and, for each triple $(m,k,c)$ of integers, the
number of flats $F$ with $|F|=m$ and $r(F)=k$ for which the
restriction $M|F$ has $c$ coloops.  Beyond the multiset of size-rank
pairs $(|A|,r(A))$ noted above as equivalent to the Tutte polynomial,
as \cite[Theorem 5.3]{cone} shows, from $\mathcal{G}(M)$, one can
find, for each triple $(m,k,c)$ of integers, the number of sets $A$
with $|A|=m$ and $r(A)=k$ for which the restriction $M|A$ has $c$
coloops.  Akin to the universality property of the Tutte polynomial
among matroid invariants that satisfy deletion-contraction rules,
Derksen and Fink \cite{valuative} showed that the
$\mathcal{G}$-invariant is a universal valuative invariant for
subdivisions of matroid base polytopes.

Reflecting on the proof of \cite[Theorem 7.3]{catdata} reveals that
the chains of cyclic flats in the configuration play the key role; the
lattice structure is secondary.  So if one can construct pairs of
non-isomorphic lattices for which one can relate their chains via
bijections, one might be able to produce pairs of matroids with
different configurations and the same $\mathcal{G}$-invariant.  That
is the idea that we develop in this paper and use to shed light on how
much stronger the configuration is compared to the
$\mathcal{G}$-invariant.  Since the Tutte polynomial can be computed
from the $\mathcal{G}$-invariant, our results also contribute to the
theory of Tutte-equivalent matroids (see \cite{survey}).

In Section \ref{sec:back}, we review the relevant background and
previously known examples of matroids with different configurations
but the same $\mathcal{G}$-invariant.  The core of the paper is
Section \ref{sec:tools}, where we develop the tools that we apply in
Sections \ref{sec:lattice} and \ref{sec:paving} to give the
constructions of interest.  These tools and the strategy in the proofs
of Theorems \ref{thm:latticeextension} and \ref {thm:generalpaving}
are likely to be useful to obtain more such constructions.

\section{Notation, background, and prior results}\label{sec:back}

Our notation and terminology for matroid theory follow Oxley
\cite{oxley}.  Let $[n]$ denote the set $\{1,2,\ldots,n\}$.  Let
$\hat{0}$ denote the least element of a lattice $L$, and let $\hat{1}$
denote its greatest element.  If we need to clarify in which lattice
an interval is formed, we use a subscript, as in $[a,b]_L$. Likewise,
we may use $\hat{0}_L$ and $\hat{1}_L$.  For a lattice $L$, we let
$L^{\circ}$ denote $L$ with $\hat{0}$ and $\hat{1}$ removed, so
$L^{\circ}$ is the open interval $(\hat{0}, \hat{1})$.

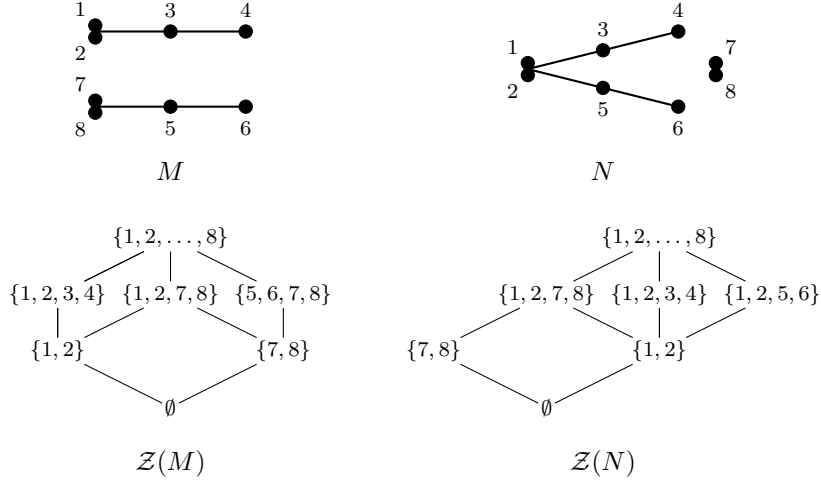
\begin{figure}
  \centering
  \begin{tikzpicture}[scale=1]
  \filldraw (0,4.92) node[below left] {\footnotesize $2$} circle  (2.5pt);
  \filldraw (0,5.08) node[above left] {\footnotesize $1$} circle  (2.5pt);
  \filldraw (1,5) node[above=2pt] {\footnotesize $3$} circle  (2.5pt);
  \filldraw (2,5) node[above=2pt] {\footnotesize $4$} circle  (2.5pt);
  \filldraw (0,3.92) node[below left] {\footnotesize $8$} circle  (2.5pt);
  \filldraw (0,4.08) node[above left] {\footnotesize $7$} circle  (2.5pt);
  \filldraw (1,4) node[below=2pt] {\footnotesize $5$} circle  (2.5pt);
  \filldraw (2,4) node[below=2pt] {\footnotesize $6$} circle  (2.5pt);
  \draw[thick](0,4)--(2,4);
  \draw[thick](0,5)--(2,5);
  \node at (1,3.15) {$M$};

  \filldraw (5.75,4.42) node[below left] {\footnotesize $2$} circle  (2.5pt);
  \filldraw (5.75,4.58) node[above left] {\footnotesize $1$} circle  (2.5pt);
  \filldraw (6.75,4.75) node[above=2pt] {\footnotesize $3$} circle  (2.5pt);
  \filldraw (7.75,5) node[above=2pt] {\footnotesize $4$} circle  (2.5pt);
  \filldraw (6.75,4.25) node[below=2pt] {\footnotesize $5$} circle  (2.5pt);
  \filldraw (7.75,4) node[below=2pt] {\footnotesize $6$} circle  (2.5pt);
  \filldraw (8.25,4.42) node[below right] {\footnotesize $8$} circle  (2.5pt);
  \filldraw (8.25,4.58) node[above right] {\footnotesize $7$} circle  (2.5pt);
  \draw[thick](5.75,4.5)--(7.75,4);
  \draw[thick](5.75,4.5)--(7.75,5);
  \node at (6.75,3.15) {$N$};

  \node[inner sep = 0.3mm] (em) at (1,0) {\footnotesize $\emptyset$};
  \node[inner sep = 0.3mm]  (s1) at (-0.5,0.75) {\footnotesize $\{1,2\}$};
  \node[inner sep = 0.3mm] (s3) at (2.5,0.75) {\footnotesize $\{7,8\}$};
  \node[inner sep = 0.3mm] (s12) at (-0.5,1.5) {\footnotesize $\{1,2,3,4\}$};
  \node[inner sep = 0.3mm] (s13) at (1,1.5) {\footnotesize $\{1,2,7,8\}$};
  \node[inner sep = 0.3mm] (s23) at (2.5,1.5) {\footnotesize $\{5,6,7,8\}$};
  \node[inner sep = 0.3mm] (s123) at (1,2.25) {\footnotesize $\{1,2,\ldots,8\}$};

  \foreach \from/\to in
  {em/s1,em/s3,s1/s12,s1/s13,s3/s23,s3/s13,s12/s123,
    s13/s123,s23/s123} \draw(\from)--(\to);

  \foreach \from/\to in{s1/s12,s12/s123} \draw (\from)--(\to);

  \node at (1,-0.75) {$\mathcal{Z}(M)$};

\node[inner sep = 0.3mm] (em) at (6,0) {\footnotesize $\emptyset$};
\node[inner sep = 0.3mm] (1) at (4.5,0.75) {\footnotesize $\{7,8\}$};
\node[inner sep = 0.3mm] (3) at (7.5,0.75) {\footnotesize $\{1,2\}$};
\node[inner sep = 0.3mm] (12) at (7.5,1.5) {\footnotesize $\{1,2,3,4\}$};
\node[inner sep = 0.3mm] (13) at (6,1.5) {\footnotesize $\{1,2,7,8\}$};
\node[inner sep = 0.3mm] (23) at (9,1.5) {\footnotesize $\{1,2,5,6\}$};
\node[inner sep = 0.3mm] (123) at (7.5,2.25) {\footnotesize $\{1,2,\ldots,8\}$};

\foreach \from/\to in
{em/1,em/3,3/12,1/13,3/23,3/13,12/123,13/123,23/123}
\draw(\from)--(\to);

\foreach \from/\to in {3/12,12/123} \draw (\from)--(\to);

  \node at (6.75,-0.75) {$\mathcal{Z}(N)$};

\end{tikzpicture}  
\caption{Two rank-$3$ matroids, $M$ and $N$, and their lattices of
  cyclic flats.}\label{fig:runningmatroid}
\end{figure}

Given a matroid $M$, a subset $A$ of $E(M)$ is \emph{cyclic} if $A$ is
a (possibly empty) union of circuits, or, equivalently, the
restriction $M|A$ has no coloops.  It follows that if $A$ and $B$ are
cyclic flats, then so are $\cl(A\cup B)$ and the flat obtained from
$A\cap B$ by removing the coloops of $M|A\cap B$; these are the join
and meet of $A$ and $B$ in the lattice $\mathcal{Z}(M)$ of cyclic
flats of $M$.  Figure \ref{fig:runningmatroid} shows two matroids and
their lattices of cyclic flats.  We will use the following elementary
result about cyclic flats.

\begin{lemma}
  Let $M$ be a matroid with neither loops nor coloops.  If $X$ is any
  cyclic flat of $M$, then $\mathcal{Z}(M|X)$ is the interval
  $[\emptyset,X]$ in $\mathcal{Z}(M)$, and
  $$\mathcal{Z}(M/X)=\{F-X\,:\,F\in \mathcal{Z}(M) \text{ and } X\subseteq
  F\},$$ so the lattice $\mathcal{Z}(M/X)$ is isomorphic to the
  interval $[X,E(M)]$ in $\mathcal{Z}(M)$.
\end{lemma}

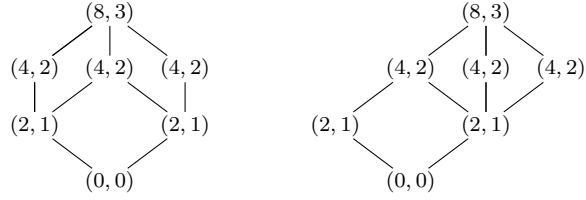
\begin{figure}
  \centering

\begin{tikzpicture}[scale=1]
  \node[inner sep = 0.3mm] (em) at (1,0) {\footnotesize $(0,0)$};
  \node[inner sep = 0.3mm]  (s1) at (0,0.75) {\footnotesize $(2,1)$};
  \node[inner sep = 0.3mm] (s3) at (2,0.75) {\footnotesize $(2,1)$};
  \node[inner sep = 0.3mm] (s12) at (0,1.5) {\footnotesize $(4,2)$};
  \node[inner sep = 0.3mm] (s13) at (1,1.5) {\footnotesize $(4,2)$};
  \node[inner sep = 0.3mm] (s23) at (2,1.5) {\footnotesize $(4,2)$};
  \node[inner sep = 0.3mm] (s123) at (1,2.25) {\footnotesize $(8,3)$};

  \foreach \from/\to in
  {em/s1,em/s3,s1/s12,s1/s13,s3/s23,s3/s13,s12/s123,
    s13/s123,s23/s123} \draw(\from)--(\to);

  \foreach \from/\to in{s1/s12,s12/s123} \draw (\from)--(\to);

\node[inner sep = 0.3mm] (em) at (5,0) {\footnotesize $(0,0)$};
\node[inner sep = 0.3mm] (1) at (4,0.75) {\footnotesize $(2,1)$};
\node[inner sep = 0.3mm] (3) at (6,0.75) {\footnotesize $(2,1)$};
\node[inner sep = 0.3mm] (12) at (6,1.5) {\footnotesize $(4,2)$};
\node[inner sep = 0.3mm] (13) at (5,1.5) {\footnotesize $(4,2)$};
\node[inner sep = 0.3mm] (23) at (7,1.5) {\footnotesize $(4,2)$};
\node[inner sep = 0.3mm] (123) at (6,2.25) {\footnotesize $(8,3)$};

\foreach \from/\to in
{em/1,em/3,3/12,1/13,3/23,3/13,12/123,13/123,23/123}
\draw(\from)--(\to);

\foreach \from/\to in {3/12,12/123} \draw (\from)--(\to);

\end{tikzpicture}  
\caption{The configurations of the matroids $M$ and $N$ in Figure
  \ref{fig:runningmatroid}.  The pairs shown are $(s(x),\rho(x))$, the
  size and rank of the corresponding cyclic flat.
}\label{fig:runningconfig}
\end{figure}

We now make the notion of the configuration, introduced informally
above, precise.  The \emph{configuration} of a matroid $M$ with no
coloops is the triple $(L, s, \rho)$, where $L$ is a lattice and
$s:L\to\mathbb{Z}$ and $\rho:L\to\mathbb{Z}$ are functions such that
there is an isomorphism $\phi: L\to \mathcal{Z}(M)$ for which
$s(x)=|\phi(x)|$ and $\rho(x)=r(\phi(x))$ for all $x\in L$. (Many
triples can satisfy these properties, but they all contain the same
data, so we view them as the same.) The configurations of the matroids
in Figure \ref{fig:runningmatroid} are shown in Figure
\ref{fig:runningconfig}. Many non-isomorphic matroids can have the
same configuration.  For instance, consider paving matroids, that is,
matroids where all flats that are properly contained in hyperplanes
are independent; any two paving matroids of the same rank, on the same
number of elements, and with the same number of dependent hyperplanes
of each size have the same configuration, and many non-isomorphic
paving matroids can share these parameters.

The configuration of $M$ gives the configurations of certain minors of
$M$, as the next lemma states.

\begin{lemma}
  Let $M$ be a matroid with neither loops nor coloops.  If
  $X\in\mathcal{Z}(M)$, then the size and rank of each set in
  $\mathcal{Z}(M|X)$ and $\mathcal{Z}(M/X)$ can be found from the
  configuration of $M$.  So, from the configuration of $M$, we get the
  configuration of each minor $M|X/Y$ with $X,Y\in \mathcal{Z}(M)$ and
  $Y\subsetneq X$.
\end{lemma}

The following result from \cite{cyclic} explains why we do not impose
any conditions on the lattices that we consider other than that they
are finite.

\begin{lemma}
  For each finite lattice $L$, there is a matroid $M$ for which
  $\mathcal{Z}(M)$ is isomorphic to $L$.
\end{lemma}

Given a matroid $M$ on $E$ with rank function $r$ and a permutation
$\pi$, say $e_1,e_2,\ldots,e_n$, of $E$, the \emph{rank sequence}
$\underline{r}(\pi)=r_1r_2\ldots r_n$ is given by
$$r_i=r(\{e_1,e_2,\ldots,e_i \})
-r(\{e_1,e_2,\ldots,e_{i-1} \}).$$ Thus, $\{e_i\,: \, r_i=1\}$ is a
basis of $M$.  Let $k=r(M)$.  Each rank sequence $\underline{r}(\pi)$
is an $(n,k)$\emph{-sequence}, that is, a sequence of $k$ ones and
$n-k$ zeroes.  For each $(n,k)$-sequence $\underline{r}$, let
$[\underline{r}]$ be a formal symbol, and let $\mathcal{G}(n,k)$ be
the vector space over a field of characteristic zero consisting of all
formal linear combinations of such symbols. The
$\mathcal{G}$-\emph{invariant of} $M$ is defined
by
$$\mathcal{G}(M)=\sum_{\text{permutations } \pi \text{ of } E}
[\underline{r}(\pi)]. $$ The matroids $M$ and $N$ in Figure
\ref{fig:runningmatroid} have the same $\mathcal{G}$-invariant, which
is
$$
384 [10101000]+
2496 [10110000]+
1344 [11001000]+
7296 [11010000]+
28800 [11100000]
$$
Among the $384$ permutations that yield the rank sequence $[10101000]$
are, for instance, $1,2,7,8,3,4,5,6$, for $M$ and $N$, and
$7,8,5,6,4,3,2,1$, for $M$ but not $N$.

In this paper we use a reformulation of $\mathcal{G}(M)$ developed in
\cite{catdata}.  An $(n,k)$\emph{-composition} is a sequence
$\boldsymbol{a}=(a_0, a_1, \ldots, a_k)$ of integers where $a_0\geq 0$
and $a_i>0$ for $i\in [k]$ such that $a_0+a_1+\cdots +a_k=n$. Let $M$
be a matroid as above.  A \emph{flag} of $M$ is a sequence
$\boldsymbol{X}=(X_0, X_1, \ldots, X_k)$ where $X_i$ is a rank-$i$
flat of $M$ and $X_i\subset X_{i+1}$ for $i<k$. The \emph{composition}
of a flag $\boldsymbol{X}$ is $\boldsymbol{a}$ where $a_0=|X_0|$ and
$a_i=|X_i-X_{i-1}|$ for $i\in [k]$. Thus, $\boldsymbol{a}$ is an
$(n,k)$-composition. Let $\nu(M;\boldsymbol{a})$ be the number of
flags of $M$ with composition $\boldsymbol{a}$. The \emph{catenary
  data} of $M$ is the $\binom{n}{k}$-dimensional vector
$(\nu(M;\boldsymbol{a}))$ that is indexed by $(n,k)$-compositions.

For example, for the matroid $M$ in Figure \ref{fig:runningmatroid},
we have
$$\nu(M;(0,1,1,6)) = 8, \quad
\nu(M;(0,1,2,5)) = 4, \quad \nu(M;(0,1,3,4)) = 4,$$
$$\nu(M;(0,2,1,5)) = 4, \quad \text{and}  \quad
\nu(M;(0,2,2,4)) = 4,
$$ 
and the same for $N$.  In $M$, each of the flats $\{1,2\}$ and
$\{7,8\}$ is in two flags that have composition $(0,2,2,4)$, while in
$N$, the flat $\{1,2\}$ is in three such flags and $\{7,8\}$ is in
only one.

A special basis of $\mathcal{G}(n,k)$, called the $\gamma$-basis, is
defined in \cite{catdata}; its vectors $\gamma(\boldsymbol{a})$ are
indexed by the $(n,k)$-compositions $\boldsymbol{a}$.  By the next
result \cite[Theorem 3.3]{catdata}, the catenary data of $M$ is
equivalent to $\mathcal{G}(M)$.

\begin{thm}\label{thm:catdata}
  The $\mathcal{G}$-invariant of $M$ is determined by its catenary
  data and conversely.  In particular,
  $$\mathcal{G}(M)=\sum_{(n,k)\text{-compositions
    }\boldsymbol{a}}\nu(M;\boldsymbol{a}) 
  \gamma(\boldsymbol{a}).$$
\end{thm}

By \cite[Theorem 7.3]{catdata}, if $M$ has no coloops, then
$\mathcal{G}(M)$ can be computed from the configuration of $M$. If $M$
has coloops, then $\mathcal{G}(M)$ can be computed from
$\mathcal{G}(M\setminus X)$ and $|X|$ where $X$ is the set of coloops.
Likewise, if $M$ has loops, then $\mathcal{G}(M)$ can be computed from
$\mathcal{G}(M\setminus Y)$ and $|Y|$ where $Y$ is the set of loops.
(More generally, if $M$ is the direct sum $M_1\oplus M_2$, then
$\mathcal{G}(M)$ can be computed from $\mathcal{G}(M_1)$ and
$\mathcal{G}(M_2)$; see \cite[Section 4.4]{catdata}.) Thus, we focus
on matroids $M$ with neither loops nor coloops, and so $\emptyset$ and
$E(M)$ are cyclic flats of $M$.  Analogous to the notation $L^{\circ}$
defined above, we let $\mathcal{Z}^{\circ}(M)$ denote the set of
nonempty, proper cyclic flats of $M$.  Thus, $\mathcal{Z}^{\circ}(M)$
is empty if and only if $M$ is a uniform matroid.

Prior to this work, the matroids with different configurations but the
same $\mathcal{G}$-invariant that we were aware of were:
\begin{itemize}
\item Dowling lattices of a given rank exceeding three, using
  non-isomorphic groups of the same order \cite{catdata},
\item the free $m$-cones of non-isomorphic matroids that have the same
  $\mathcal{G}$-invariant, along with tipless, baseless, and
  tipless/baseless variations \cite{cone}, and
\item the duals of any such examples.
\end{itemize}
(The effect of taking the dual of a matroid on the
$\mathcal{G}$-invariant, in the original formulation, is simple:
reverse the order within each rank sequence, and then switch $0$s and
$1$s .)

Another construction is easy to treat using catenary data: given a
matroid $M$ and a positive integer $t$, construct the extension $M_t$
of $M$ by, for each $e \in E(M)$, adding a set $A_e$ of $t$ elements
parallel to $e$, with $A_e$ disjoint from all other such sets and from
$E(M)$.  (Contracting the tip of a free $m$-cone of $M$ yields $M_m$.)
Note that all flats of $M_t$ are cyclic, $\mathcal{Z}(M_t)$ is
isomorphic to the lattice of flats of $M$, and each flag
$(Y_0, Y_1, \ldots, Y_k)$ of $M_t$ corresponds to a flag
$(X_0, X_1, \ldots, X_k)$ of $M$ where $X_i\subseteq Y_i$ and
$|Y_i|= |X_i|(t+1)$ for each $i$.  So if $M$ and $N$ are
non-isomorphic and have the same $\mathcal{G}$-invariant, then $M_t$
and $N_t$ have the same catenary data, and so the same
$\mathcal{G}$-invariant, but different configurations.

For any positive integer $n$, there is a positive integer $m$ so that
there are at least $n$ non-isomorphic groups of order $m$, so Dowling
lattices yield sets of arbitrarily many matroids with different
configurations (specifically, non-isomorphic lattices of cyclic flats)
but the same $\mathcal{G}$-invariant.  Our constructions yield other
large sets of matroids of this type for which the lattices of cyclic
flats are simpler and smaller (see Examples \ref{ex:1} and
\ref{ex:2}).

\section{The tools}\label{sec:tools}

The first lemma generalizes an argument from the proof of
\cite[Theorem 7.3]{catdata}.  The join $F_1\vee F_2$ mentioned below
is the join in the lattice $\mathcal{Z}(M)$ of cyclic flats.

\begin{lemma}\label{lem:simplify}
  Assume that the function $g:\mathcal{Z}^{\circ}(M)\to 2^{2^{E(M)}}$ has the
  property that if $g(F_1)\cap g(F_2)\ne \emptyset$, then
  $F_1\vee F_2\in\mathcal{Z}^{\circ}(M)$ and
  $g(F_1)\cap g(F_2) \subseteq g(F_1\vee F_2)$.  Then
  $$\Bigl|\bigcup_{F\in \mathcal{Z}^{\circ}(M)}g(F) \Bigr|  =
  \sum_{\substack{S\subseteq \mathcal{Z}^{\circ}(M), \\
      S\ne\emptyset}} (-1)^{|S|+1} \Bigl|\bigcap_{F\in S}g(F) \Bigr| =
  \sum_{\substack{\text{\emph{nonempty chains }} \\
      S\subseteq \mathcal{Z}^{\circ}(M)}} (-1)^{|S|+1}\Bigl|\bigcap_{F\in S}g(F)
  \Bigr|.$$
\end{lemma}

Before proving the lemma, we note an example of such a function $g$.
For a matroid $M$, let $g_M:\mathcal{Z}^{\circ}(M)\to 2^{2^{E(M)}}$ be
given by
\begin{equation}\label{eq:gdef}
  g_M(F) = \{X\,:\,X\subset E(M),\, |X|=r(M)-1, \text{ and }\cl_M(X\cap
  F)=F\}.
\end{equation}
For example, for the matroid $M$ in Figure \ref{fig:runningmatroid},
the sets in $g_M(\{1,2,3,4\})$ are precisely the five bases of
$M|\{1,2,3,4\}$, while those in $g_M(\{1,2\})$ are the thirteen
$2$-element sets that contain at least one of $1$ and $2$.  It is
routine to verify that, for any matroid $M$, the function $g_M$ has
the property in the hypothesis of the lemma.  Later in this section we
make heavy use of the function $g_M$ for various matroids $M$.

\begin{proof}[Proof of Lemma \ref{lem:simplify}]
  The principle of inclusion/exclusion gives the equality of the first
  two expressions, so we focus on the two sums.  Let $F_1$ and $F_2$
  be incomparable flats in $\mathcal{Z}^{\circ}(M)$ with
  $g(F_1)\cap g(F_2)\ne \emptyset$, so $F_1\vee F_2$ properly contains
  both of them.  Now
  $g(F_1)\cap g(F_2) =g(F_1)\cap g(F_2)\cap g(F_1\vee F_2)$ by the
  assumed property.  If $F_1$ and $F_2$ are in a subset $S$ of
  $\mathcal{Z}^{\circ}(M)$, and if $S'$ is the symmetric difference
  $S\triangle\{F_1\vee F_2\}$, then $(-1)^{|S|}=-(-1)^{|S'|}$ and
  $$\bigcap_{F\in S}g(F) = \bigcap_{F\in S'}g(F),$$ so such terms
  in the first sum could cancel.  To cancel all such terms, take a
  linear extension $\leq$ of the order $\subseteq$ on
  $\mathcal{Z}^{\circ}(M)$ and, if a subset $S$ of
  $\mathcal{Z}^{\circ}(M)$ contains incomparable cyclic flats, let
  $F_1$ and $F_2$ be such a pair for which $(F_1,F_2)$ is least in the
  lexicographic order that $\leq$ induces on
  $\mathcal{Z}^{\circ}(M)\times\mathcal{Z}^{\circ}(M)$; in the first
  sum, cancel the term that arises from $S$ with the one that arises
  from $S\triangle\{F_1\vee F_2\}$.  Any incomparable cyclic flats
  $F_1$ and $F_2$ in a set $S$ that remains after cancellation satisfy
  $g(F_1)\cap g(F_2)= \emptyset$; such terms contribute zero to the
  first sum and so can be omitted.  The result of all such
  cancellations and omissions is the sum on the right side.
\end{proof}

Let $(X_0, X_1, \ldots, X_r)$ be a flag of a matroid $M$ of rank $r$
that has neither loops nor coloops, so $X_0=\emptyset$ and $X_r=E(M)$.
If, for each flat $X_i$, we remove the coloops of $M|X_i$ from $X_i$,
we obtain a chain of cyclic flats.  The (possibly empty) chain
obtained by removing $\emptyset$ and $E(M)$ from this chain is the
\emph{reduced cyclic chain} of the flag.  For a chain $C$ in
$\mathcal{Z}^\circ(M)$, let $\fl(C)$ be the set of all flags of
$M$ whose reduced cyclic chain is $C$.  For a set $T$ of chains in
$\mathcal{Z}^\circ(M)$, let $\fl(T)$ be the union of all sets $\fl(C)$
where $C\in T$.  For a set $S$ of flags, let $\comp(S)$ be the
multiset of compositions of the flags in $S$.  Thus, $\comp(\fl(T))$
is the multiset of all compositions of all of the flags of $M$ whose
reduced cyclic chain is in the set $T$ of chains in
$\mathcal{Z}^\circ(M)$.

\begin{example}\label{ex:original}
  For the reduced cyclic chain $C=(\{1,2\})$ in the matroid $M$ in
  Figure \ref{fig:runningmatroid}, $\fl_M(C)$ contains four flags,
  namely, (omitting $\emptyset$ and $E(M)$ for brevity)
  $(\{1,2\}, \{1,2,i\})$ and $(\{i\}, \{1,2,i\})$ for $i\in \{5,6\}$;
  the compositions in $\comp(\fl_M(C))$ are $(0,2,1,5)$ and
  $(0,1,2,5)$, and each has multiplicity $2$. For the matroid $N$ in
  that figure, $\fl_N(C)$ is empty.  Likewise, for $C'=(\{7,8\})$,
  $\fl_M(C')$ contains four flags: $(\{7,8\}, \{7,8,i\})$ and
  $(\{i\}, \{7,8,i\})$ for $i\in \{3,4\}$; the compositions in
  $\comp(\fl_M(C'))$ are $(0,2,1,5)$ and $(0,1,2,5)$, each of which
  has multiplicity $2$. Also, $\fl_N(C')$ contains eight flags:
  $(\{7,8\}, \{7,8,i\})$ and $(\{i\}, \{7,8,i\})$ for
  $i\in \{3,4,5,6\}$; the compositions in the multiset
  $\comp(\fl_N(C'))$ are $(0,2,1,5)$ and $(0,1,2,5)$, and each has
  multiplicity $4$.  Thus, the multisets $\comp(\fl_M(C))$ and
  $\comp(\fl_M(C'))$ differ from $\comp(\fl_N(C))$ and
  $\comp(\fl_N(C'))$, but
  $\comp(\fl_M(\{C,C'\}))=\comp(\fl_N(\{C,C'\}))$, so together $C$ and
  $C'$ make the same contribution to the catenary data of $M$ and $N$.
  
  To complete that example, note that neither $M$ nor $N$ has flags
  with reduced cyclic chain $(\{1,2,7,8\})$.  For each of $M$ and $N$,
  there are eight flags with the empty reduced cyclic chain, and each
  has composition $(0,1,1,6)$.  In each of $M$ and $N$, each of the
  four $2$-element reduced cyclic chains arises from exactly one flag,
  and that flag has composition $(0,2,2,4)$.  Finally, the reduced
  cyclic chains $(\{1,2,3,4\})$ and $(\{5,6,7,8\})$ in $M$, and
  $(\{1,2,3,4\})$ and $(\{1,2,5,6\})$ in $N$, each arise from two
  flags (e.g., $(\{3\}, \{1,2,3,4\})$ and $(\{4\}, \{1,2,3,4\})$ for
  the first); each of these reduced cyclic chains contributes two
  copies of $(0,1,3,4)$ to the catenary data.
\end{example}

This example motivates a general technique, stated next, for showing
that two matroids have the same $\mathcal{G}$-invariant.  This result,
which follows immediately from Theorem \ref{thm:catdata} and the
definitions above, is what we use in Sections \ref{sec:lattice} and
\ref{sec:paving} to show that the constructions presented there
produce matroids with the same $\mathcal{G}$-invariant.

\begin{thm}\label{thm:technique}
  Let $M$ and $N$ be matroids.  Let $\{P_1,P_2,\ldots,P_d\}$ be a
  partition of the set of reduced cyclic chains of
  $\mathcal{Z}^{\circ}(M)$. Let $\{Q_1,Q_2,\ldots,Q_d\}$ be a
  partition of the set of reduced cyclic chains of
  $\mathcal{Z}^{\circ}(N)$.  If
  $\comp(\fl_M(P_i)) =\comp(\fl_N(Q_i))$ for each $i\in[d]$, then
  $\mathcal{G}(M) = \mathcal{G}(N)$.
\end{thm}

In the rest of this section, we treat the tools that we will use to
establish equalities of the form
$\comp(\fl_{M}(P_i)) =\comp(\fl_{N}(Q_i))$.

Part of the proof of \cite[Theorem 7.3]{catdata} establishes the
following lemma.  In this result and what follows, we use $\iota(N)$
to denote the number of independent hyperplanes (that is, hyperplanes
that are independent sets) of a matroid $N$.

\begin{lemma}\label{lem:flags}
  Let $M$ be a matroid that has no loops and no coloops, so
  $\emptyset$ and $E(M)$ are cyclic flats of $M$.  Let
  $C=\{F_1\subsetneq F_2\subsetneq \cdots\subsetneq F_t\}$ be a chain
  in $\mathcal{Z}^{\circ}(M)$.  Set $F_0=\emptyset$ and
  $F_{t+1}=E(M)$.  Let $\mathcal{L}$ be the set of all lists of length
  $r(M)$, with no repeated entries, that are obtained as follows:
  \begin{itemize}
  \item for each $j\in[t+1]$, pick an independent hyperplane $H_j$ of
    $M|F_j/F_{j-1}$,
  \item consider lists in which the entries are the sets
    $F_1,F_2,\ldots,F_{t+1}$ along with the singleton subsets of
    $H_1,H_2,\ldots,H_{t+1}$,
  \item such a list is in $\mathcal{L}$ if and only if, for each
    $j\in[t]$, $F_j$ occurs before $F_{j+1}$, and, for each
    $j\in[t+1]$, all singleton subsets of $H_j$ occur before $F_j$ .
  \end{itemize}
  Mapping $L \in \mathcal{L}$ to the sequence $\phi(L)$, the $i$th
  entry of which is the union of the first $i$ sets in $L$, for
  $0\leq i \leq r(M)$, defines a bijection
  $\phi:\mathcal{L}\to \fl(C)$.
  
  Thus, the multiset $\comp(\fl(C))$ is determined by the sizes and
  ranks of $F_1,F_2,\ldots,F_t$, and $E(M)$, along with
  $\iota(M|F_1)$, $\iota(M|F_j/F_{j-1})$ for $2\leq j\leq t$, and
  $\iota(M/F_t)$.
\end{lemma}

The number of independent hyperplanes of a matroid can be
deduced from its $\mathcal{G}$-invariant \cite[Proposition
5.5]{catdata} and so from its configuration.

\begin{cor}\label{cor:chainsimplify}
  Let $M$ and $M'$ be matroids, both having neither loops nor coloops.
  Let $C= \{F_1\subsetneq F_2\subsetneq \cdots\subsetneq F_t\}$ be a
  chain in $\mathcal{Z}^{\circ}(M)$, and 
  $C'=\{F'_1\subsetneq F'_2\subsetneq \cdots\subsetneq F'_t\}$ a chain
  in $\mathcal{Z}^{\circ}(M')$. Assume that the following pairs of
  minors have the same configuration:
  \begin{itemize}
  \item[(a)] the restrictions $M|F_1$ and $M'|F'_1$,
  \item[(b)] for each $i$ with $2\leq i\leq t$, the minors $M|F_i/F_{i-1}$
    and $M'|F'_i/F'_{i-1}$, and
  \item[(c)] the contractions $M/F_t$ and $M'/F'_t$.
  \end{itemize}
  Then $\comp(\fl_M(C))=\comp(\fl_{M'}(C'))$.
\end{cor}

Corollary \ref{cor:chainsimplify} and Lemmas \ref{lem:configsofchains}
and \ref{lem:grouping} are the key tools for the results in the rest
of the paper.  We next prove a result that we use in the proofs of
those lemmas.

\begin{lemma}\label{lem:gtwochains}
  Let $M$ and $M'$ be matroids of rank $r$ for which $|E(M)|=|E(M')|$.
  Assume that $M$ and $M'$ have neither loops nor coloops.  Let
  $\{F_1\subsetneq F_2\subsetneq \cdots\subsetneq F_t\}$ be a chain in
  $\mathcal{Z}^\circ(M)$, and
  $\{F'_1\subsetneq F'_2\subsetneq \cdots\subsetneq F'_t\}$ a chain in
  $\mathcal{Z}^\circ(M')$, for which the following pairs of minors
  have the same configuration:
  \begin{itemize}
  \item[(a)] the restrictions $M|F_1$ and $M'|F'_1$, and
  \item[(b)] for each $i$ with $2\leq i\leq t$, the minors $M|F_i/F_{i-1}$
    and $M'|F'_i/F'_{i-1}$.
  \end{itemize}
  Let $g_M$ and $g_{M'}$ be given by equation \emph{(\ref{eq:gdef})}.
  Then
  \begin{equation}\label{eq:gtwochains}
    \Bigl|\bigcap_{i\in[t]}g_M(F_i) \Bigr| =
    \Bigl|\bigcap_{i\in[t]}g_{M'}(F'_i) \Bigr|.
  \end{equation}
\end{lemma}

\begin{proof}
  We get each set in $g_M(F_1)\cap g_M(F_2)\cap\cdots\cap g_M(F_t)$
  exactly once by following the steps below, allowing for all possible
  choices:
  \begin{itemize}
  \item choose integers $a_1,a_2,\ldots,a_t,a_{t+1}$ with
    $a_1+a_2+\cdots+a_t+a_{t+1}=r-1$, and with $a_1\geq r_M(F_1)$ and
    $a_i\geq r(F_i)-r(F_{i-1})$ for $i$ with $2\leq i\leq t$, and
    $a_{t+1}\geq 0$;
  \item choose
    \begin{itemize}
    \item[(i)] a spanning set $W_1$ of $M|F_1$ where $|W_1|=a_1$,
    \item [(ii)] for each $i$ with $2\leq i\leq t$, a spanning set
      $W_i$ of $M|F_i/F_{i-1}$ where $|W_i|=a_i$, and
    \item [(iii)] any set $W_{t+1}\subseteq E(M)-F_t$ where
      $|W_{t+1}|=a_{t+1}$;
    \end{itemize}
  \item then
    $W_1\cup W_2\cup\cdots\cup W_{t+1}\in g_M(F_1)\cap
    g_M(F_2)\cap\cdots\cap g_M(F_t)$.
  \end{itemize}
  A similar sequences of choices yields all sets in
  $g_{M'}(F'_1)\cap g_{M'}(F'_2)\cap\cdots\cap g_{M'}(F'_t)$.  From
  this, equality (\ref{eq:gtwochains}) follows since, by the
  hypothesis, the restrictions $M|F_1$ and $M'|F'_1$ have the same
  configuration and so the same number of spanning sets of size $a_1$,
  and, for each $i$ with $2\leq i\leq t$, the minors $M|F_i/F_{i-1}$
  and $M'|F'_i/F'_{i-1}$ have the same configuration and so the same
  number of spanning sets of size $a_i$, and those assumptions along
  with the hypothesis $|E(M)|=|E(M')|$ also give
  $|E(M)-F_t|=|E(M')-F'_t|$.
\end{proof}

For a matroid $M$, let $\ch(\mathcal{Z}^{\circ}(M))$ be the set of
non-empty chains in $\mathcal{Z}^{\circ}(M)$.

The next two lemmas and their proofs have a similar flavor and could
be merged.  We treat them separately for clarity.  Note that for the
matroids $M$ and $N$ in Figure \ref{fig:runningmatroid}, there are
many bijections
$\sigma: \ch(\mathcal{Z}^{\circ}(M))\to \ch(\mathcal{Z}^{\circ}(N))$
that satisfy the conditions below, but none is induced by a bijection
from $\mathcal{Z}^{\circ}(M)$ onto $\mathcal{Z}^{\circ}(N)$.

\begin{lemma}\label{lem:configsofchains}
  Let the rank-$r$ matroids $M$ and $M'$ on $n$ elements have neither
  loops nor coloops.  Assume that there is a bijection
  $\sigma: \ch(\mathcal{Z}^{\circ}(M))\to
  \ch(\mathcal{Z}^{\circ}(M'))$ such that for each chain
  $S=\{X_1\subsetneq X_2\subsetneq \cdots\subsetneq X_t\}$ in
  $\ch(\mathcal{Z}^{\circ}(M))$,
  \begin{itemize}
  \item[(a)] its image
    $\sigma(S)=\{X'_1\subsetneq X'_2\subsetneq \cdots\subsetneq X'_t\}$
    has $|S|$ elements, and
  \item[(b)] the restrictions $M|X_1$ and $M'|X'_1$ have the same
    configuration, as do, for each $i$ with $2\leq i\leq t$, the
    minors $M|X_i/X_{i-1}$ and $M'|X'_i/X'_{i-1}$.
  \end{itemize}
  Then $\iota(M)=\iota(M')$ and, for the empty chain $\emptyset$,
  $\comp(\fl_M(\emptyset)) = \comp(\fl_{M'}(\emptyset))$.
\end{lemma}

\begin{proof}
  Let $g_M$ and $g_{M'}$ be given by equation (\ref{eq:gdef}).  Then
  $$\iota(M)=\binom{n}{r-1}-\Bigl|\bigcup_{F\in
    \mathcal{Z}^{\circ}(M)} g_M(F) \Bigr|,$$ and similarly for $M'$,
  so in order to show that $\iota(M)=\iota(M')$, it suffices to show
  that
  $$\Bigl|\bigcup_{F\in \mathcal{Z}^{\circ}(M)} g_M(F) \Bigr|
  =\Bigl|\bigcup_{F\in \mathcal{Z}^{\circ}(M')} g_{M'}(F) \Bigr|.$$ By
  Lemma \ref{lem:simplify}, it suffices to show that
  $$\sum_{S\in  \ch(\mathcal{Z}^{\circ}(M))} (-1)^{|S|+1}\Bigl|\bigcap_{F\in
    S}g_M(F) \Bigr| = \sum_{S\in \ch(\mathcal{Z}^{\circ}(M'))}
  (-1)^{|S|+1}\Bigl|\bigcap_{F\in S}g_{M'}(F) \Bigr|.$$ This holds by
  applying Lemma \ref{lem:gtwochains} to each pair of chains
  $(S,\sigma(S))$ for $S\in \ch(\mathcal{Z}^{\circ}(M))$. 

  We now get $\comp(\fl_M(\emptyset)) = \comp(\fl_{M'}(\emptyset))$
  since the reduced cyclic chain of a flag is empty if and only if the
  composition of the flag is $(0,1,1,\ldots,1,n-r+1)$, and, by Lemma
  \ref{lem:flags}, the number of such flags in $M$ is
  $(r-1)!\iota(M)$, and likewise for $M'$.
\end{proof}

In Example \ref{ex:original}, we have
$\comp(\fl_M(\{C,C'\}))=\comp(\fl_N(\{C,C'\}))$ for reduced cyclic
chains $C$ and $C'$ even though what $C$ and $C'$ contribute
individually differs in $M$ versus $N$.  The next lemma identifies
conditions that yield the same conclusion.

\begin{lemma}\label{lem:grouping}
  Let $M$ and $M'$ be rank-$r$ matroids on $n$ elements with no loops
  and no coloops.  Let $C_1,C_2,\ldots,C_p$ be distinct (not
  necessarily disjoint) nonempty chains in $\mathcal{Z}^{\circ}(M)$
  and let $C'_1,C'_2,\ldots,C'_p$ be distinct (not necessarily
  disjoint) nonempty chains in $\mathcal{Z}^{\circ}(M')$ with
  $|C_i|=|C'_j|$ for all $i,j\in[p]$.  Write $C_i$ as
  $\{F_{i,1}\subsetneq F_{i,2}\subsetneq\cdots\subsetneq F_{i,t}\}$
  and $C'_i$ as
  $\{F'_{i,1}\subsetneq F'_{i,2}\subsetneq\cdots\subsetneq
  F'_{i,t}\}$.  Assume that for all $i,j\in[p]$,
  \begin{itemize}
  \item[(a)] $M|F_{i,1}$ and $M'|F'_{j,1}$ have the same
    configuration, as do the minors $M|F_{i,h}/F_{i,h-1}$ and
    $M'|F'_{j,h}/F'_{j,h-1}$ for each $h$ with $2\leq h\leq t$.
  \end{itemize}
  Also assume that there is a bijection
  $$\Phi:\bigcup_{i\in[p]}\bigl(\ch(\mathcal{Z}^{\circ}(M/F_{i,t}))\times\{C_i\}
  \bigr) \to \bigcup_{i\in[p]}
  \bigl(\ch(\mathcal{Z}^{\circ}(M'/F'_{i,t}))\times\{C'_i\}\bigr)$$
  with the following property:
  \begin{itemize}
  \item[(b)] if
    $\Phi( (\{X_1\subsetneq X_2\subsetneq\cdots\subsetneq X_k\},C_i))
    = (\{X'_1\subsetneq X'_2\subsetneq\cdots\subsetneq
    X'_{k'}\},C'_j)$, then $k=k'$ and $M|X_1/F_{i,t}$ and
    $M'|X'_1/F'_{j,t}$ have the same configuration, as do the minors
    $M|X_h/X_{h-1}$ and $M'|X'_h/X'_{h-1}$ for each $h$ with
    $2\leq h \leq k$.
  \end{itemize}
  Then
  $\comp(\fl_M(\{C_1,C_2,\ldots,C_p\})) =
  \comp(\fl_{M'}(\{C'_1,C'_2,\ldots,C'_p\}))$.
\end{lemma}

This result implies Corollary \ref{cor:chainsimplify}, but it applies
more broadly, even for $p=1$, since property (b) may hold for some
bijection from $\ch(\mathcal{Z}^{\circ}(M/F_{1,t}))$ onto
$\ch(\mathcal{Z}^{\circ}(M'/F'_{1,t}))$ even if $M/F_{1,t}$ and
$M'/F'_{1,t}$ have different configurations (see the proof of Theorem
\ref{thm:latticeextension}).  We may have $F_{i,t} = F_{j,t}$ even if
$i\ne j$, so chains in $\ch(\mathcal{Z}^{\circ}(M/F_{i,t}))$ may arise
multiple times; the second entry in the Cartesian product
distinguishes these occurrences.  When $p=1$, we can omit the second
factor in the Cartesian product.

\begin{proof}[Proof of Lemma \ref{lem:grouping}]
  By Lemma \ref{lem:flags}, we get each flag of $M$ whose reduced
  cyclic chain is $C_i$ exactly once by the following procedure:
  \begin{itemize}
  \item pick an independent hyperplane of $M|F_{i,1}$ and take a
    permutation of its singleton subsets; insert the set $F_{i,1}$ as
    the last entry of this list;
  \item treat each $h$ with $2\leq h\leq t$ in order; for $h$, pick
    an independent hyperplane of $M|F_{i,h}/F_{i,h-1}$ and insert its
    singleton subsets anywhere into the list we have so far, and then
    insert the set $F_{i,h}$ as the last entry of this list;
  \item pick an independent hyperplane of $M/F_{i,t}$ and insert its
    singleton subsets anywhere into the list we have so far, and then
    insert $E(M)$ as the last entry of the list;
  \end{itemize}
  the flag that we get from the final list of sets has as its $j$th flat, for
  $0\leq j\leq r$, the union of the first $j$ sets in the list.  The
  conditions above make it clear that the multiset of compositions of
  flags that result from any particular choice of independent
  hyperplanes of $M|F_{i,1}$, $M|F_{i,h}/F_{i,h-1}$ for
  $2\leq h\leq t$, and $M/F_{i,t}$ does not depend on those sets; only
  their sizes (which are fixed) and the sizes of
  $F_{i,1}, F_{i,2},\ldots, F_{i,t}$ and $E(M)$ matter.  A similar
  description applies to the flags of $M'$ whose reduced cyclic chain
  is $C'_j$.  The numbers of independent hyperplanes of $M|F_{i,1}$
  and of $M|F_{i,h}/F_{i,h-1}$ for $2\leq h\leq t$ are determined by
  their configurations, so by condition (a) they are equal to the
  numbers of independent hyperplanes of the corresponding minors of
  $M'$ using the chain $C'_j$.  The hypotheses do not imply that
  $\iota(M/F_{i,t})$ and $\iota(M'/F_{j,t})$ are equal, but for the
  conclusion in the theorem to hold, all that we need is
  \begin{equation}\label{eq:sumiota}
    \sum_{i\in[p]}\iota(M/F_{i,t})=    \sum_{i\in[p]}\iota(M'/F'_{i,t}).
  \end{equation}

  For each $i\in[p]$, let $g_{M/F_{i,t}}$, which we shorten to $g_i$,
  be given by equation (\ref{eq:gdef}), and likewise for $g'_i$, the
  shortened form of $g_{M'/F'_{i,t}}$.  Then
  $$\iota(M/F_{i,t})=\binom{n-|F_{i,t}|}{r-r_M(F_{i,t})-1}-\,
  \Bigl|\bigcup_{F\in \mathcal{Z}^{\circ}(M/F_{i,t})} g_i(F) \Bigr|,$$
  and similarly for $\iota(M'/F'_{j,t})$.  We get
  $|F_{i,t}| = |F'_{j,t}|$ and $r_M(F_{i,t})= r_{M'}(F'_{j,t})$ for
  all $i,j\in[p]$ by condition (a), so equation (\ref{eq:sumiota})
  will follow by showing that
  $$\sum_{i\in[p]}\,\, \Bigl|\bigcup_{F\in \mathcal{Z}^{\circ}(M/F_{i,t})} g_i(F)
  \Bigr| = \sum_{i\in[p]}\,\,\Bigl|\bigcup_{F\in
    \mathcal{Z}^{\circ}(M'/F'_{i,t})} g'_i(F) \Bigr|.$$ By Lemma
  \ref{lem:simplify}, it suffices to show that
  $$\sum_{i\in[p]} \sum_{S\in \ch(\mathcal{Z}^{\circ}(M/F_{i,t}))}
  (-1)^{|S|+1}\Bigl|\bigcap_{F\in S}g_i(F) \Bigr| = \sum_{i\in[p]}
  \sum_{S\in \ch(\mathcal{Z}^{\circ}(M'/F'_{i,t}))}
  (-1)^{|S|+1}\Bigl|\bigcap_{F\in S}g'_i(F) \Bigr|.$$ Lemma
  \ref{lem:gtwochains} applies to the chains $S$ and $S'$ whenever
  $\Phi((S,C_i)) = (S',C'_j)$ with $i,j\in[p]$, so this equality
  holds.
\end{proof}

\section{Application: a construction modifying a
  lattice}\label{sec:lattice}

In this section we show how to extend a finite lattice in different
ways to produce lattices in different configurations that yield the
same $\mathcal{G}$-invariant.

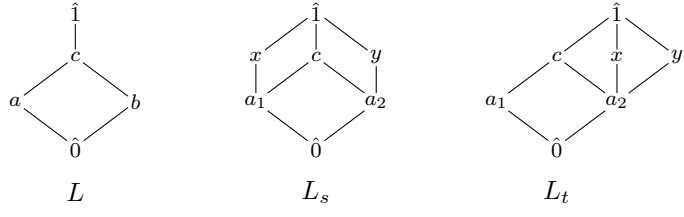
\begin{figure}
  \centering
  \begin{tikzpicture}[scale=0.8]

  \node[inner sep = 0.3mm] (0em) at (-3,0) {\footnotesize $\hat{0}$};
  \node[inner sep = 0.3mm]  (0s1) at (-4,0.75) {\footnotesize $a$};
  \node[inner sep = 0.3mm] (0s3) at (-2,0.75) {\footnotesize $b$};
  \node[inner sep = 0.3mm] (0s13) at (-3,1.5) {\footnotesize $c$};
  \node[inner sep = 0.3mm] (0s123) at (-3,2.25) {\footnotesize $\hat{1}$};

  \foreach \from/\to in {0em/0s1,0em/0s3,0s1/0s13,0s3/0s13,0s13/0s123}
  \draw(\from)--(\to);

  \node at (-3,-0.75) {$L$};
  
  \node[inner sep = 0.3mm] (em) at (1,0) {\footnotesize $\hat{0}$};
  \node[inner sep = 0.3mm]  (s1) at (0,0.75) {\footnotesize $a_1$};
  \node[inner sep = 0.3mm] (s3) at (2,0.75) {\footnotesize $a_2$};
  \node[inner sep = 0.3mm] (s12) at (0,1.5) {\footnotesize $x$};
  \node[inner sep = 0.3mm] (s13) at (1,1.5) {\footnotesize $c$};
  \node[inner sep = 0.3mm] (s23) at (2,1.5) {\footnotesize $y$};
  \node[inner sep = 0.3mm] (s123) at (1,2.25) {\footnotesize $\hat{1}$};

  \foreach \from/\to in
  {em/s1,em/s3,s1/s12,s1/s13,s3/s23,s3/s13,s12/s123,
    s13/s123,s23/s123} \draw(\from)--(\to);

  \node at (1,-0.75) {$L_s$};
  
  \node[inner sep = 0.3mm] (sem) at (5,0) {\footnotesize $\hat{0}$};
  \node[inner sep = 0.3mm] (1) at (4,0.75) {\footnotesize $a_1$};
  \node[inner sep = 0.3mm] (3) at (6,0.75) {\footnotesize $a_2$};
  \node[inner sep = 0.3mm] (12) at (6,1.5) {\footnotesize $x$};
  \node[inner sep = 0.3mm] (13) at (5,1.5) {\footnotesize $c$};
  \node[inner sep = 0.3mm] (23) at (7,1.5) {\footnotesize $y$};
  \node[inner sep = 0.3mm] (123) at (6,2.25) {\footnotesize $\hat{1}$};

  \foreach \from/\to in
  {sem/1,sem/3,3/12,1/13,3/23,3/13,12/123,13/123,23/123}
  \draw(\from)--(\to);

    \node at (5,-0.75) {$L_t$};
\end{tikzpicture}
\caption{A lattice $L$ and two extension of the type in the
  construction treated in Section \ref{sec:lattice}, where $L_s$ is
  isomorphic to $\mathcal{Z}(M)$ and where $L_t$ is isomorphic to
  $\mathcal{Z}(N)$ for the matroids $M$ and $N$ in Figure
  \ref{fig:runningmatroid}.  }\label{fig:firstex}
\end{figure}

We first construct the lattices that appear in Theorem
\ref{thm:latticeextension}.  Let $a_1,a_2,\ldots,a_m$ be distinct
elements of a finite lattice $L$ for which
\begin{itemize}
\item there is a $b\in L$ with $a_i\wedge a_j=b$ for all $i,j\in[m]$
  with $i\ne j$, and
\item for each $i\in [m-1]$, there is a lattice isomorphism
  $\tau_{m,i}:[\hat{0},a_m]\to [\hat{0},a_i]$ with $\tau_{m,i}(y)=y$
  for all $y\in [\hat{0},b]$.
\end{itemize}
Set $\tau_{i,m} = \tau_{m,i}^{-1}$ and, for distinct $i,j\in[m-1]$,
set $\tau_{i,j} = \tau_{m,j}\circ\tau_{i,m}$.  Therefore
$\tau_{i,j}:[\hat{0},a_i]\to [\hat{0},a_j]$ is a lattice isomorphism
and $\tau_{i,j}(y)=y$ for all $y\in [\hat{0},b]$.  Let
$L_1,L_2,\ldots,L_n$ be finite lattices that are disjoint from each
other and from $L$.  Fix two functions $s:[n]\to [m]$ and
$t:[n]\to [m]$.  Form a lattice $L_s$ as follows:
\begin{itemize}
\item for $i\in[n]$, form $L'_i$ from $L_i$ by replacing
  $\hat{1}_{L_i}$ by $\hat{1}_L$ and $\hat{0}_{L_i}$ by $a_{s(i)}$,
\item viewing lattices as relations, let $L_s$ be the transitive
  closure of $L\cup L'_1\cup L'_2\cup\cdots\cup L'_n$.
\end{itemize}
It is routine to check that $L_s$ is indeed a lattice.  In terms of
Hasse diagrams, $L_s$ is obtained by, for each $i\in [n]$, inserting
$L_i$ into the interval $[a_{s(i)},\hat{1}_L]$ of $L$, where
$\hat{0}_{L_i}$ is identified with $a_{s(i)}$, and $\hat{1}_{L_i}$ is
identified with $\hat{1}_L$.  Define $L_t$ similarly. Thus, $L_s$ and
$L_t$ have the same elements.  The lattices $L_s$ and $L_t$ may be
isomorphic, but in many examples they are not.  Figure
\ref{fig:firstex} shows how the lattices of cyclic flats of the two
matroids in Figure \ref{fig:runningmatroid} are obtained by extending
a lattice using this construction.  We give more examples after the
proof of the next result.

\begin{thm}\label{thm:latticeextension}
  Let $L_s$ and $L_t$ be as defined above.  Let $M_s$ and $M_t$ be
  matroids, with rank functions $r_s$ and $r_t$, respectively, neither
  having loops nor coloops, for which, for some lattice isomorphisms
  $\phi_s:L_s \to \mathcal{Z}(M_s)$ and
  $\phi_t:L_t \to \mathcal{Z}(M_t)$,
  \begin{itemize}
  \item[(1)] $|\phi_s(y)|=|\phi_t(y)|$ and
    $r_s(\phi_s(y))=r_t(\phi_t(y))$ for all $y$ in $L_s$, and
  \item[(2)] if $y\in [\hat{0},a_m]$, then
    $|\phi_s(y)|=|\phi_s(\tau_{m,i}(y))|$ and
    $r_s(\phi_s(y))=r_s(\phi_s(\tau_{m,i}(y)))$ for all $i\in [m-1]$.
  \end{itemize}
  Then $\mathcal{G}(M_s)=\mathcal{G}(M_t)$.
\end{thm}

From the $\mathcal{G}$-invariant of a matroid, we can deduce the
number of cyclic flats of each size and rank \cite[Proposition
5.5]{catdata}, so condition (1) just says how this data lines up.
Condition (2) is the only restrictive assumption on this data.
 
\begin{proof}[Proof of Theorem \ref{thm:latticeextension}]
  By Theorem \ref{thm:technique}, it suffices produce a partition
  $\{P_1,P_2,\ldots,P_d\}$ of the set of chains in
  $\mathcal{Z}^{\circ}(M_s)$ and a partition $\{Q_1,Q_2,\ldots,Q_d\}$
  of the set of chains in $\mathcal{Z}^{\circ}(M_t)$ so that, for each
  $h\in[d]$, one of Corollary \ref{cor:chainsimplify}, Lemma
  \ref{lem:configsofchains}, or Lemma \ref{lem:grouping} gives
  $\comp(\fl_{M_s}(P_h)) =\comp(\fl_{M_t}(Q_h))$.  We specify chains
  in $\mathcal{Z}^\circ(M_s)$ as chains in $L^{\circ}_s$; to get the
  corresponding chains in $\mathcal{Z}^\circ(M_s)$, apply $\phi_s$; we
  handle chains in $\mathcal{Z}^\circ(M_t)$ similarly.  Let
  $\ch(L^{\circ}_s)$ denote the set of nonempty chains in
  $L^{\circ}_s$; define $\ch(L^{\circ}_t)$ similarly.
  
  We start with some useful observations for verifying that certain
  minors have the same configuration, as required by Corollary
  \ref{cor:chainsimplify} and Lemmas \ref{lem:configsofchains} and
  \ref{lem:grouping}.  These observations are cast in terms of an
  interval $[x,z]$ in $L_s$.  The corresponding minor
  $M_s|\phi_s(z)/\phi_s(x)$ of $M_s$ can be, if $x=\hat{0}$, the
  restriction to the first cyclic flat in a chain, or, if $z=\hat{1}$,
  the contraction by the last cyclic flat in a chain, or the minor
  determined by consecutive cyclic flats in a chain.  Only
  applications of Corollary \ref{cor:chainsimplify} use intervals of
  the form $[x,\hat{1}]$.  Properties (1) and (2) imply that, for all
  distinct integers $i,j\in [m]$ and all $y\in [\hat{0},a_i]$,
  \begin{itemize}
  \item[(S)]
    $|\phi_s(y)|=|\phi_s(\tau_{i,j}(y))|=|\phi_t(\tau_{i,j}(y))|$,
    and
  \item[(R)]
    $r_s(\phi_s(y))=r_s(\phi_s(\tau_{i,j}(y)))=
    r_t(\phi_t(\tau_{i,j}(y)))$.
  \end{itemize}
  Consider an interval $[x,z]_{L_s}$ with $x\ne z$.  The next three
  items treat all intervals that arise below when we apply Corollary
  \ref{cor:chainsimplify} and Lemmas \ref{lem:configsofchains} and 
  \ref{lem:grouping}.
  \begin{itemize}
  \item If there is no $i\in[m]$ with $x\leq a_i<z$, then
    $M_s|\phi_s(z)/\phi_s(x)$ and $M_t|\phi_t(z)/\phi_t(x)$ have the
    same configuration.
  \end{itemize}
  This holds since $[x,z]_{L_s}$ is also an interval of $L_t$ and the
  equalities in property (1) hold for all $y$ in $[x,z]_{L_s}$.  This
  next statement gives addition information for some of the intervals
  to which the last statement applies.
  \begin{itemize}
  \item If $z\leq a_i$ for some $i\in[m]$, then, for all $j\in[m]$
    with $j\ne i$, the minors $M_s|\phi_s(z)/\phi_s(x)$ and
    $M_t|\phi_t(\tau_{i,j}(z))/\phi_t(\tau_{i,j}(x))$ have the same
    configuration.
  \end{itemize}
  This holds since properties (S) and (R) apply to all
  $y\in [x,z]$.
  \begin{itemize}
  \item If $x\leq a_i$ and $z\in L^\circ_h$ where $s(h)=i$ and
    $t(h)=j$, then
    \begin{itemize}
    \item[$\circ$] if $i=j$, then the minors $M_s|\phi_s(z)/\phi_s(x)$
      and $M_t|\phi_t(z)/\phi_t(x)$ have the same configuration, while
    \item[$\circ$] if $i\ne j$, then the minors
      $M_s|\phi_s(z)/\phi_s(x)$ and
      $M_t|\phi_t(z)/\phi_t(\tau_{i,j}(x))$ have the same
      configuration.
    \end{itemize}
  \end{itemize}
  This holds by property (1) and, for the second part, properties (S)
  and (R) applied to all $y\in [x,a_i]$.
  
  To treat the empty chain via Lemma \ref{lem:configsofchains},
  consider $\sigma: \ch(L^{\circ}_s)\to \ch(L^{\circ}_t)$ given by,
  for $C \in \ch(L^{\circ}_s)$,
  $$\sigma(C) =  \begin{cases}
    (C-L)\cup\tau_{s(h),t(h)}(C\cap L), & \text{if }
    C\cap L_h^{\circ}\ne\emptyset \text{ and } s(h)\ne t(h),\\
    C, & \text{otherwise.}
  \end{cases}$$ It is routine to check that $\sigma$ is a bijection.
  Condition (a) clearly holds.  The observations in the second
  paragraph show that condition (b) holds, so Lemma
  \ref{lem:configsofchains} applies.

  A chain $C$ is in the symmetric difference
  $\ch(L^{\circ}_s)\triangle \ch(L^{\circ}_t)$ if and only if (i)
  $C\not\subseteq L$, (ii) $C\cap L\not\subseteq[\hat{0},b]$, and
  (iii) $s(h)\ne t(h)$ for the $h$ with $C-L\subseteq L^{\circ}_h$.
  Restricting $\sigma$ gives a bijection
  $\sigma: \ch(L^{\circ}_s)-\ch(L^{\circ}_t)\to
  \ch(L^{\circ}_t)-\ch(L^{\circ}_s)$.  Corollary
  \ref{cor:chainsimplify} applies to the chains $\phi_s(C)$ and
  $\phi_t(\sigma(C))$ for each
  $C\in \ch(L^{\circ}_s)-\ch(L^{\circ}_t)$; the observations in the
  second paragraph show that conditions (a)--(c) hold.

  We now focus on chains in $\ch(L^{\circ}_s)\cap \ch(L^{\circ}_t)$.
  Let $\max(C)$ denote the greatest element in a chain $C$.

  Let $C$ be a nonempty chain of both $L^{\circ}_s$ and $L^{\circ}_t$.
  If no $i\in[m]$ has $C\subseteq [\hat{0},a_i]$, then Corollary
  \ref{cor:chainsimplify} applies to $\phi_s(C)$ and $\phi_t(C)$, with
  the observations in the second paragraph showing that conditions
  (a)--(c) hold.  If $C\subseteq [\hat{0},b]$, then Lemma
  \ref{lem:grouping} applies to $\phi_s(C)$ and $\phi_t(C)$ (so
  $p=1$): restrict the map $\sigma$ defined above to
  $ \ch((\max(C),\hat{1})_{L^{\circ}_s})$ to get a bijection
  $\Phi: \ch((\max(C),\hat{1})_{L^{\circ}_s})\to
  \ch((\max(C),\hat{1})_{L^{\circ}_t})$; conditions (a) and (b) hold
  by the observations in the second paragraph.  Now consider a chain
  $C$ for which there is exactly one $i\in [m]$ with
  $C\subseteq [\hat{0},a_i]$.  The images of $C$ under the maps
  $\tau_{i,j}$, along with $C$ itself, give $m$ chains
  $C_1,C_2,\ldots,C_m$ with $C_j\subseteq [\hat{0},a_j]$ and
  $\tau_{j,h}(C_j)=C_h$ for all distinct $j,h\in[m]$.  We apply Lemma
  \ref{lem:grouping} to the chains
  $\phi_s(C_1),\phi_s(C_2),\ldots,\phi_s(C_m)$ in
  $\mathcal{Z}^{\circ}(M_s)$, and
  $\phi_t(C_1),\phi_t(C_2),\ldots,\phi_t(C_m)$ in
  $\mathcal{Z}^{\circ}(M_t)$.  By the observations in the second
  paragraph, condition (a) in that lemma holds.  Note that if
  $D\in \ch((\max(C_h),\hat{1})_{L_s})$ and
  $D\cap L^\circ_j\ne\emptyset$, then $s(j)=h$.  The map
  $$\Phi:\bigcup_{h\in[m]}
  \bigl(\ch((\max(C_h),\hat{1})_{L_s}) \times\{C_h\}\bigr) \to
  \bigcup_{h\in[m]} \bigl(\ch((\max(C_h),\hat{1})_{L_t})
  \times\{C_h\}\bigr)$$ given by
  $$\Phi((D,C_h)) = \begin{cases}
    ((D-L)\cup\tau_{h,t(j)}(D\cap L) ,C_{t(j)}), & \text{if }
    D\cap L_j^{\circ}\ne\emptyset \text{ and } t(j)\ne h,\\
    (D,C_h), & \text{otherwise}, \end{cases}$$ is easily seen to be a
  bijection that, by the observations in the second paragraph, when we
  apply $\phi_s$ and $\phi_t$, satisfies property (b) in Lemma
  \ref{lem:grouping}.
\end{proof}

\begin{figure}
    \centering
    \begin{tikzpicture}[scale=1]

  \node[inner sep = 0.3mm] (e) at (3.75,0) {\footnotesize $\hat{0}$};
  \node[inner sep = 0.3mm] (a1) at (0.75,0.75) {\footnotesize $a_1$};
  \node[inner sep = 0.3mm] (a2) at (3.25,0.75) {\footnotesize $a_2$};
  \node[inner sep = 0.3mm] at (5,0.75) {\footnotesize $\cdots$};
  \node[inner sep = 0.3mm] (am) at (6.75,0.75) {\footnotesize $a_m$};

  \node[inner sep = 0.3mm] (x1) at (0,1.5) {\footnotesize $x_1$};
  \node[inner sep = 0.3mm] at (0.75,1.5) {\footnotesize $\cdots$};
  \node[inner sep = 0.3mm] (x2) at (1.5,1.5) {\footnotesize $x_{i-1}$};

  \node[inner sep = 0.3mm] (x3) at (2.5,1.5) {\footnotesize $x_i$};
  \node[inner sep = 0.3mm] at (3.25,1.5) {\footnotesize $\cdots$};
  \node[inner sep = 0.4mm] (x4) at (4,1.5) {\footnotesize $x_{j-1}$};
  \node[inner sep = 0.3mm] at (5,1.5) {\footnotesize $\cdots$};
  
  \node[inner sep = 0.3mm] (x5) at (6,1.5) {\footnotesize $x_k$};
  \node[inner sep = 0.3mm] at (6.75,1.5) {\footnotesize $\cdots$};
  \node[inner sep = 0.4mm] (x6) at (7.5,1.5) {\footnotesize $x_n$};

  \node[inner sep = 0.3mm] (t) at (3.75,2.25) {\footnotesize $\hat{1}$};
  \foreach \from/\to in
  {e/a1,e/a2,e/am,a1/x1,a1/x2,a2/x3,a2/x4,am/x5,am/x6,t/x1,t/x2,t/x3,t/x4,t/x5,t/x6}
  \draw(\from)--(\to);
    \end{tikzpicture}
  \caption{The lattice in Example \ref{ex:1}.}\label{fig:ex1}
\end{figure}
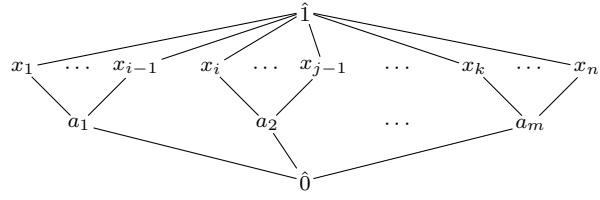

\begin{example}\label{ex:1}
  If the lattice $L$ in the construction above is isomorphic to the
  lattice of flats of an $m$-point line and each $L_i$ is a
  three-element chain, then $L_s$ has the form illustrated in Figure
  \ref{fig:ex1}, where some open intervals $(a_h,\hat{1})$ may be
  empty.  The number of non-isomorphic lattices of this type is the
  number of integer partitions of $n$ with at most $m$ parts.  One
  type of matroid for which its lattice of cyclic flats has this form
  is a rank-$4$ matroid with $m$ three-point lines
  $A_1,A_2,\ldots,A_m$, each pair of which spans the matroid, and
  where each cyclic plane contains just one cyclic line, namely, some
  $A_i$.  If the sizes of the $n$ cyclic planes are fixed and
  distinct, then, up to isomorphism, $\sum_{i\in[m]}S(n,i)$ matroids
  satisfy these conditions, and their configurations are different.
  (Here, $S(n,i)$ is the Stirling number of the second kind.  For some
  pairs of configurations, the lattices are isomorphic but the
  assignments of sizes differ.)  This produces huge sets of matroids
  with very simple but different configurations and the same
  $\mathcal{G}$-invariant.  For instance, for $m=n$, the cyclic planes
  can have distinct sizes with only $n(n+9)/2$ elements in the
  matroids, and the size of the set of matroids produced this way,
  with different configurations and the same $\mathcal{G}$-invariant,
  is the $n$th Bell number.  More such examples result by altering the
  sizes and ranks.
\end{example}

\begin{example}\label{ex:2}
  Since a cyclic flat is both cyclic (a union of circuits) and a flat
  (an intersection of hyperplanes, the complements of which are
  circuits of the dual matroid), it follows that $X$ is a cyclic flat
  of a matroid $M$ if and only if $E(M)-X$ is a cyclic flat of the
  dual matroid, $M^*$.  Thus, $\mathcal{Z}(M^*)$ is isomorphic to the
  order dual of $\mathcal{Z}(M)$.  So a counterpart of Theorem
  \ref{thm:latticeextension} holds for the order duals of $L_s$ and
  $L_t$.  An example of a matroid whose lattice of cyclic flats is
  isomorphic to the order dual of the lattice considered in the
  previous example is a rank-$r$ matroid, with $r\geq 4$, in which the
  restriction to each of the cyclic hyperplanes $A_1,A_2,\ldots,A_m$,
  which are pairwise disjoint and have the same cardinality, is a
  paving matroid, where any two cyclic flats of rank $r-2$ in distinct
  cyclic hyperplanes span the matroid.  It is conjectured that
  asymptotically almost all matroids are paving, and, while that is
  still unresolved (see \cite{asymp} for recent results), the number
  of paving matroids is known to be enormous, so this construction
  produces huge sets of matroids with different configurations but the
  same $\mathcal{G}$-invariant.  Also, as in the previous example, as
  $r$ grows, we have more options for the ranks of the cyclic flats
  $A_i$ and for the cyclic flats $X_j$ that they contain.
\end{example}

\section{Application: a construction using paving
  matroids}\label{sec:paving}

In Theorem \ref{thm:generalpaving} we show how to use the lattices of
flats of paving matroids as the lattices in different configurations
that yield the same $\mathcal{G}$-invariant.  Examples are given after
the proof.

\begin{thm}\label{thm:generalpaving}
  Let $N_1$ and $N_2$ be non-isomorphic rank-$r$ paving matroids on a
  set $E$.  For $i\in[2]$, let $\mathcal{H}_i$ be the set of
  hyperplanes of $N_i$ and let $L_i$ be the lattice of flats of $N_i$.
  Assume that there are partitions $\{A_1,A_2,\ldots,A_p\}$ of
  $\mathcal{H}_1- \mathcal{H}_2$ and $\{B_1,B_2,\ldots,B_p\}$ of
  $\mathcal{H}_2- \mathcal{H}_1$, and, for each $j\in [p]$, a
  bijection $\alpha_j:E\to E$ such that $X\in A_j$ if and only if
  $\alpha_j(X)\in B_j$.  Let $M_1$ and $M_2$ be matroids, with rank
  functions $r_1$ and $r_2$, respectively, for which there are lattice
  isomorphisms $\phi_i:L_i\to \mathcal{Z}(M_i)$, for $i\in[2]$, such
  that
  \begin{itemize}
  \item if $X\in L_1\cap L_2$, then $|\phi_1(X)|=|\phi_2(X)|$ and
    $r_1(\phi_1(X))=r_2(\phi_2(X))$, and
  \item if $j\in[p]$ and $Y\in A_j$, then
    \begin{itemize}
    \item[$\circ$] $|\phi_1(Y)|=|\phi_2(\alpha_j(Y))|$ and
      $r_1(\phi_1(Y))=r_2(\phi_2(\alpha_j(Y)))$, and
    \item[$\circ$] $|\phi_1(X)|=|\phi_1(\alpha_j(X))|$ and
      $r_1(\phi_1(X))=r_1(\phi_1(\alpha_j(X)))$ for all $X\subseteq Y$
      with $|X|<r-1$.
    \end{itemize}
  \end{itemize}
  Then $\mathcal{G}(M_1)=\mathcal{G}(M_2)$.  
\end{thm}

The hypotheses imply that $N_1$ and $N_2$ have the same configuration.
Also, the conditions on $\phi_1$ and $\phi_2$ give the following
equalities: if $Y\in A_j$ and $X\subseteq Y$ with $|X|<r-1$, then
\begin{itemize}
\item
  $|\phi_2(X)| =|\phi_1(X)|=|\phi_1(\alpha_j(X))|=
  |\phi_2(\alpha_j(X))|$ and
\item
  $r_2(\phi_2(X)) =r_1(\phi_1(X)) =r_1(\phi_1(\alpha_j(X)))
  =r_2(\phi_2(\alpha_j(X)))$.
  \end{itemize}

\begin{proof}[Proof of Theorem \ref{thm:generalpaving}]
  We use the technique outlined in Theorem \ref{thm:technique}: we
  give partitions $\{P_1,P_2,\ldots,P_d\}$ of the set of chains in
  $\mathcal{Z}^{\circ}(M_1)$ and $\{Q_1,Q_2,\ldots,Q_d\}$ of the set
  of chains in $\mathcal{Z}^{\circ}(M_2)$ so that, for each $h\in[d]$,
  we get $\comp(\fl_{M_1}(P_h)) =\comp(\fl_{M_2}(Q_h))$ from Corollary
  \ref{cor:chainsimplify}, Lemma \ref{lem:configsofchains}, or Lemma
  \ref{lem:grouping}.  We can specify chains in
  $\mathcal{Z}^{\circ}(M_i)$, for $i\in[2]$, as chains in
  $L^{\circ}_i$; to get the corresponding chain in
  $\mathcal{Z}^{\circ}(M_i)$, apply $\phi_i$.  Let $\ch(L^{\circ}_i)$
  denote the set of nonempty chains in $L^{\circ}_i$.

  To treat the empty chain via Lemma \ref{lem:configsofchains}, we
  need a bijection $\sigma: \ch(L^{\circ}_1)\to \ch(L^{\circ}_2)$ so
  that the composition $\phi_2\circ\sigma\circ\phi_1^{-1}$ satisfies
  the conditions in Lemma \ref{lem:configsofchains}.  For a chain
  $C=\{Y_1\subsetneq Y_2\subsetneq \cdots \subsetneq Y_t\}$ in
  $\ch(L^{\circ}_1)$, let
  $$\sigma(C) =  \begin{cases}
    \alpha_j(C), & \text{if } Y_t\in A_j,\\
    C, & \text{otherwise,}
  \end{cases}$$ where $\alpha_j(C) = \{\alpha_j(Y_h)\,:\,h\in[t]\}$.
  The hypotheses imply that $\sigma$ is a bijection.  Condition (a)
  clearly holds.  To see that $\phi_2\circ\sigma\circ\phi_1^{-1}$
  satisfies condition (b), set $Y_0=\emptyset$.  If $\sigma(C)=C$,
  then $|\phi_1(X)| =|\phi_2(X)|$ and $r_1(\phi_1(X)) =r_2(\phi_2(X))$
  for all $X\in[Y_{h-1},Y_h]$ with $h\in[t]$, which gives condition
  (b) in this case.  If $\sigma(C)=\alpha_j(C)$, then
  $|\phi_1(X)| =|\phi_2(\alpha_j(X))|$ and
  $r_1(\phi_1(X)) =r_2(\phi_2(\alpha_j(X)))$ for all
  $X\in[Y_{h-1},Y_h]$ with $h\in[t]$, which gives condition (b) in
  this case.
  
  Now
  $\ch(L^{\circ}_{i_1})-\ch(L^{\circ}_{i_2}) = \bigl\{\bigr
  \{Y_1\subsetneq Y_2\subsetneq \cdots \subsetneq Y_t\} \,:\, Y_t\in
  \mathcal{H}_{i_1}-\mathcal{H}_{i_2} \}$ for $(i_1,i_2)$ in
  $\{(1,2),(2,1)\}$.  For each
  $C = \{Y_1\subsetneq Y_2\subsetneq \cdots \subsetneq Y_t\} \in
  \ch(L^{\circ}_1)-\ch(L^{\circ}_2)$, there is a unique $j\in [p]$
  with $Y_t\in A_j$, and Corollary \ref{cor:chainsimplify} applies to
  the chains $\phi_1(C)$ of $M_1$ and $\phi_2(\alpha_j(C))$ of $M_2$.
  Properties (a)--(c) hold as in the previous paragraph.

  We now treat the nonempty chains in
  $\ch(L^{\circ}_1)\cap \ch(L^{\circ}_2)$.  For any chain $C$ that
  contains a set in $\mathcal{H}_1\cap\mathcal{H}_2$, Corollary
  \ref{cor:chainsimplify} applies to the chains $\phi_1(C)$ in $M_1$
  and $\phi_2(C)$ in $M_2$.  To treat the remaining chains, we define
  a graph $G$ whose vertices are the nonempty chains in
  $\ch(L^{\circ}_1)\cap \ch(L^{\circ}_2)$ that contain no hyperplane.
  Two such chains
  $C=\{Y_1\subsetneq Y_2\subsetneq \cdots \subsetneq Y_t\}$ and
  $C'=\{Y'_1\subsetneq Y'_2\subsetneq \cdots \subsetneq Y'_{t'}\}$ are
  adjacent in $G$ if $t=t'$ and there is a $j\in [p]$ and an
  $H\in A_j$ such that either (i) $Y_t\subseteq H$ and
  $C'=\alpha_j(C)$ or (ii) $Y'_t\subseteq H$ and $C=\alpha_j(C')$.
  Thus, if $C$ and $C'$ are adjacent, then, setting
  $Y_0=Y'_0= \emptyset$, there is a bijection
  $ \beta_h:[Y_{h-1},Y_h]\to[Y'_{h-1},Y'_h] $ for each $h\in [t]$ such
  that for each $X\in[Y_{h-1},Y_h]$,
  \begin{itemize}
  \item
    $|\phi_2(X)| =|\phi_1(X)|=|\phi_1(\beta_h(X))|= |\phi_2(\beta_h(X))|$, and
  \item
    $r_2(\phi_2(X)) =r_1(\phi_1(X)) =r_1(\phi_1(\beta_h(X)))
    =r_2(\phi_2(\beta_h(X)))$.
  \end{itemize}
  (The bijection $\beta_h$ is a restriction of $\alpha_j$ if
  $Y_t\subseteq H$ and $C'=\alpha_j(C)$, otherwise $\beta_h$ is a
  restriction of $\alpha^{-1}_j$.)  So the following minors of $M_1$
  and $M_2$ have the same configuration:
  \begin{itemize}
  \item the restrictions $M_2|\phi_2(Y_1)$, $M_1|\phi_1(Y_1)$,
    $M_1|\phi_1(Y'_1)$, and $M_2|\phi_2(Y'_1)$, and
  \item the minors $M_2|\phi_2(Y_h)/\phi_2(Y_{h-1})$,
    $M_1|\phi_1(Y_h)/\phi_1(Y_{h-1})$,
    $M_1|\phi_1(Y'_h)/\phi_1(Y'_{h-1})$, and
    $M_2|\phi_2(Y'_h)/\phi_2(Y'_{h-1})$ for each $h$ with
    $2\leq h\leq t$.
  \end{itemize}
  It follows that the same conclusions hold for any two chains in the
  same component of $G$.  Let $C_1,C_2,\ldots,C_s$ be the vertices in
  a connected component of $G$.  We complete the proof by showing that
  Lemma \ref{lem:grouping} applies to
  $\phi_1(C_1),\phi_1(C_2),\ldots,\phi_1(C_s)$ in $M_1$ and
  $\phi_2(C_1),\phi_2(C_2),\ldots,\phi_2(C_s)$ in $M_2$.  Property (a)
  in that lemma follows from what we just deduced.  For $k\in [s]$,
  let $Z_k$ be the largest set in the chain $C_k$.  The sets
  $Z_1,Z_2,\ldots, Z_s$ need not be distinct.  Consider the map
  $$\Phi:\bigcup_{k\in [s]}\left(\ch((Z_k,E)_{L_1})\times\{C_k\}
  \right)\to \bigcup_{k\in
    [s]}\left(\ch((Z_k,E)_{L_2})\times\{C_k\}\right)$$ given as
  follows: for a chain
  $D=\{W_1 \subsetneq W_2 \subsetneq\cdots \subsetneq W_s\}$ in the
  interval $(Z_k,E)_{L_1}$,
  $$\Phi((D,C_k)) =
  \begin{cases} (\alpha_j(D),\alpha_j(C_k)), & \text{if } W_s\in
    A_j,\\
    (D,C_k), & \text{otherwise}.
  \end{cases}$$ It is easy to see that $\Phi$ is a bijection.  When we
  apply $\phi_1$ and $\phi_2$, property (b) in Lemma
  \ref{lem:grouping} follows by the same type of argument as in the
  second paragraph.
\end{proof}

When two paving matroids $N_1$ and $N_2$ on $E$ have the same
configuration, partitions $\{A_1,A_2,\ldots,A_p\}$ and
$\{B_1,B_2,\ldots,B_p\}$ and bijections $\alpha_j:E\to E$, for
$j\in[p]$, that satisfy the hypothesis of Theorem
\ref{thm:generalpaving} exist since the blocks can be singletons and
the multiset of cardinalities of the hyperplanes will be the same for
$N_1$ as for $N_2$.  However, poorly-chosen partitions and bijections
can restrict the size-rank data in $\mathcal{Z}(M_1)$ and
$\mathcal{Z}(M_2)$ more than necessary.  We illustrate two efficient
choices of partitions and bijections, as well as how the matroids
$M_1$ and $M_2$ can be constructed.

\begin{example}\label{ex:3}
  Fix integers $m\geq 2$ and $n\geq 2$.  Let $N_1$ and $N_2$ be the
  rank-$3$ paving matroids on the set
  $\{a,b,x_1,x_2,\ldots,x_m,y_1,y_2,\ldots,y_n\}$ where 
  \begin{itemize}
  \item the dependent hyperplanes of $N_1$ are
    $\{a,x_1,x_2,\ldots,x_m\}$ and $\{b,y_1,y_2,\ldots,y_n\}$, and
  \item the dependent hyperplanes of $N_2$ are
    $\{a,x_1,x_2,\ldots,x_m\}$ and $\{a,y_1,y_2,\ldots,y_n\}$.
  \end{itemize} 
  These paving matroids are shown in Figure \ref{fig:ex3}.  Now
  \begin{itemize}
  \item
    $\mathcal{H}_1-\mathcal{H}_2 = \bigl\{ \{b,y_1,y_2,\ldots,y_n\},
    \{a,y_1\}, \{a,y_2\},\ldots , \{a,y_n\} \bigr\} $, and
  \item
    $\mathcal{H}_2-\mathcal{H}_1 = \bigl\{ \{a,y_1,y_2,\ldots,y_n\},
    \{b,y_1\}, \{b,y_2\},\ldots , \{b,y_n\} \bigr\}.$
  \end{itemize}
  We partition each of these into a single block.  A bijection
  $\alpha$ that has the properties in Theorem \ref{thm:generalpaving}
  is the $2$-cycle $(a,b)$.

  \begin{figure}
  \centering
  \begin{tikzpicture}[scale=0.95]
  
  \filldraw (0,1) node[above]{\footnotesize$a$} circle  (2.5pt);
  \filldraw (1,1) node[above] {\footnotesize$x_1$} circle  (2.5pt);
  \filldraw (4,1) node[above] {\footnotesize$x_m$} circle  (2.5pt);
  \filldraw (3,1) node[above] {\footnotesize$x_{m-1}$} circle  (2.5pt);
  \filldraw (0,0) node[above] {\footnotesize$b$} circle  (2.5pt);
  \filldraw (1,0) node[above] {\footnotesize$y_1$} circle  (2.5pt);
  \filldraw (4,0) node[above] {\footnotesize$y_n$} circle  (2.5pt);
  \filldraw (3,0) node[above] {\footnotesize$y_{n-1}$} circle  (2.5pt);
  \draw[thick](0,0)--(1.5,0);
  \draw[thick](0,1)--(1.5,1);
  \node at (2,0) {$\ldots$};
  \node at (2,1) {$\ldots$};
  \draw[thick](2.5,0)--(4,0);
  \draw[thick](2.5,1)--(4,1);
  
  \node at (2,-0.5) {$N_1$};

  \filldraw (6,0.5) node[above] {\footnotesize{$a$}} circle  (2.5pt);
  \filldraw (7,0.65) node[above] {\footnotesize{$x_1$}} circle  (2.5pt);
  \filldraw (9,0.85) node[above] {\footnotesize{$x_{m-1}$}} circle  (2.5pt);
  \filldraw (10,1) node[above] {\footnotesize{$x_m$}} circle  (2.5pt);
  \filldraw (7,0.35) node[below] {\footnotesize{$y_1$}} circle  (2.5pt);
  \filldraw (9,0.15) node[below] {\footnotesize{$y_{n-1}$}} circle  (2.5pt);
  \filldraw (10,0) node[below] {\footnotesize{$y_n$}} circle  (2.5pt);
  \filldraw (10.5,0.5) node[right]{\footnotesize{$b$}} circle  (2.5pt);
  
  \draw[thick](6,0.5)--(7.5,0.3);
  \draw[thick](6,0.5)--(7.5,0.7);

  \node[rotate=10] at (8,0.75) {$\ldots$};
  \node[rotate=-10] at (8,0.25) {$\ldots$};
  
  \draw[thick](8.5,0.2)--(10,0);
  \draw[thick](8.5,0.8)--(10,1);

  \node at (8,-0.5) {$N_2$};  
  \end{tikzpicture}
  \caption{The rank-$3$ paving matroids in Example
    \ref{ex:3}.}\label{fig:ex3}
  \end{figure}  

  With this pair of rank-$3$ paving matroids $N_1$ and $N_2$, we
  illustrate how straightforward it is to construct matroids $M_1$
  and $M_2$ that satisfy the hypotheses of Theorem
  \ref{thm:generalpaving} and so have different configurations and the
  same $\mathcal{G}$-invariant.  It is well known that a matroid $M$
  is determined by the set $\{(X,r(X))\,:\,X\in \mathcal{Z}(M)\}$.  We
  construct $M_1$ and $M_2$ by giving their cyclic flats and the ranks
  of these sets; this approach is justified by the following result
  \cite[Theorem 3.2]{cyclic}; see also \cite{sims}.

  \begin{lemma}
    For a set $\mathcal{Z}$ of subsets of a set $E$ and a function
    $r : \mathcal{Z}\to \mathbb{Z}$, there is a matroid $M$ on $E$
    with $\mathcal{Z}(M) = \mathcal{Z}$ and $r_M(X) = r(X)$ for all
    $X\in\mathcal{Z}$ if and only if
    \begin{itemize}
    \item[(Z0)] $(\mathcal{Z},\subseteq)$ is a lattice,
    \item[(Z1)] $r(0_\mathcal{Z})=0$, where $0_\mathcal{Z}$ is the
      least set in $\mathcal{Z}$,
    \item[(Z2)] $0<r(Y)-r(X)<|Y-X|$ for all $X,Y \in \mathcal{Z}$ with
      $X\subsetneq Y$, and
    \item[(Z3)]
      $r(X\vee Y)+r(X\wedge Y)+|(X\cap Y)-(X\wedge Y)|\leq r(X)+r(Y)$
      for $X, Y \in \mathcal{Z}$.
    \end{itemize}
  \end{lemma}

  Let $S=\{A,B,X_1,X_2,\ldots,X_m,Y_1,Y_2,\ldots,Y_n\}$ be a set of
  pairwise disjoint sets, each with at least seven elements, and with
  $|A|=|B|$.  Let $a\in E(N_1)$ correspond to $A\in S$, and likewise
  for the other elements and sets.  Let $E$ be the union of the sets
  in $S$.  Let $\mathcal{Z}_1$ be the set of unions of sets in $S$
  that correspond to flats of $N_1$, and let $\mathcal{Z}_2$ be the
  set of unions of sets in $S$ that correspond to flats of
  $N_2$. Thus, $\mathcal{Z}_1$ and $\mathcal{Z}_2$ are lattices.
  Define $r_1:\mathcal{Z}_1\to\mathbb{Z}$ and
  $r_2:\mathcal{Z}_2\to\mathbb{Z}$ as follows.
  \begin{itemize}
  \item Set $r_1(\emptyset)=r_2(\emptyset)=0$.
  \item Set $r_1(E)=r_2(E)=10$.
  \item For each set $Z\in S$, pick an $i\in\{5,6\}$ and set
    $r_1(Z)=r_2(Z)=i$, where, in addition,
    $r_1(A)=r_2(A)=r_1(B)=r_2(B)$.
  \item Pick an $i\in\{8,9\}$ and set $r_1(A\cup B)=r_2(A\cup B)=i$.
  \item For each $j\in[m]$, pick an $i\in\{8,9\}$ and set
    $r_1(B\cup X_j)=r_2(B\cup X_j)=i$.
  \item For each $j\in[n]$, pick an $i\in\{8,9\}$ and set
    $r_1(A\cup Y_j)=r_2(B\cup Y_j)=i$.
  \item For each $j\in[m]$ and $k\in [n]$, pick an $i\in\{8,9\}$ and set
    $r_1(X_j\cup Y_k)=r_2(X_j\cup Y_k)=i$.
  \item Pick an $i\in\{8,9\}$ and set
    $r_1(A\cup X_1\cup\cdots\cup X_m )=r_2(A\cup X_1\cup\cdots\cup
    X_m)=i$.
  \item Pick an $i\in\{8,9\}$ and set
    $r_1(B\cup Y_1\cup\cdots\cup Y_n )=r_2(A\cup Y_1\cup\cdots\cup
    Y_n)=i$.
  \end{itemize}
  It is easy to see that properties (Z0)--(Z3) hold for the pairs
  $(\mathcal{Z}_1,r_1)$ and $(\mathcal{Z}_2,r_2)$, so this yields
  matroids $M_1$ and $M_2$ to which Theorem \ref{thm:generalpaving}
  applies.

  Even within the narrow numerical range in this example, there can be
  more flexibility than presented above.  For instance,
  \begin{itemize}
  \item if a set in $S$ has rank five, it might have only six
    elements;
  \item a set in $S$ can have rank seven (in which case it must have
    at least eight elements) provided that it is not the intersection
    of two cyclic flats that each have rank eight;
  \item sets in $S$ could have rank four provided that, for any two
    such sets, sets that contain their union and correspond to lines
    of $N_1$ or $N_2$ have rank eight;
  \item if proper care is taken when assigning ranks (so that
    condition (Z3) will hold), the sets in $S$ do not have to be
    pairwise disjoint (for a smaller example of this type, compare
    $N_1$ and $N_2$);
  \item while we considered just certain unions of the sets in $S$, we
    could allow certain supersets of these unions.
  \end{itemize}
  Combining these factors with allowing other ranks yields an
  abundance of matroid $M_1$ and $M_2$ to which Theorem
  \ref{thm:generalpaving} applies, and the same can be done for any
  paving matroid $N_1$ and $N_2$.
\end{example}

In the concluding example, we focus exclusively on an efficient choice
of the partitions $\{A_1,A_2,\ldots,A_p\}$ and
$\{B_1,B_2,\ldots,B_p\}$ and the bijections $\alpha_j:E\to E$, for
$j\in[p]$.

\begin{figure}
  \centering
\begin{tikzpicture}[scale=1]
\draw[thick, black!20]
(0,0)--(3,0)--(3.6,0.6)--(0.6,0.6)--(0.4,2.3)--(0,0);
\filldraw[white] (3,0.5) rectangle  (3.2,0.7);
\draw[thick, black!20](0.6,0.6)--(0,0);
\draw[thick, black!20](3,0)--(3.4,2.3)--(3.6,0.6);
\draw[thick, black!20](3.4,2.3)--(0.4,2.3);

\filldraw (0.4,2.3) node[above] {\footnotesize $a$} circle  (2.2pt);
\filldraw (3.4,2.3) node[above] {\footnotesize $b$} circle  (2.2pt);
\filldraw (0.6,0.6) node[left] {\footnotesize $c$} circle  (2.2pt);
\filldraw (3.6,0.6) node[right] {\footnotesize $d$} circle  (2.2pt);
\filldraw (0,0) node[below] {\footnotesize $e$} circle  (2.2pt);
\filldraw (3,0) node[below] {\footnotesize $f$} circle  (2.2pt);
\filldraw (-0.6,0.2) node[left] {\footnotesize $p$} circle  (2.2pt);
\filldraw (-0.3,1.2) node[left] {\footnotesize $q$} circle  (2.2pt);
\filldraw (-0.1,2.1) node[left] {\footnotesize $r$} circle  (2.2pt);
\filldraw (-0.3,1.6) node[left] {\footnotesize $s$} circle  (2.2pt);
\filldraw (-0.7,0.7) node[left] {\footnotesize $t$} circle  (2.2pt);
\filldraw (-0.8,1.8) node[left] {\footnotesize $u$} circle  (2.2pt);

  \node at (1.75,-0.75) {$N_1$};  
\end{tikzpicture}
\hspace{40pt}
\begin{tikzpicture}[scale=1]
  \draw[thick, black!20] (0,0)--(3,0)--(3.6,0.6)--(0.6,0.6)--(0,0);
  \draw[thick, black!20] (0.5,1)--(3.5,1)--(4.1,1.6)--(1.1,1.6)--(0.5,1);
  \draw[thick, black!20] (0,2)--(3,2)--(3.6,2.6)--(0.6,2.6)--(0,2);

\filldraw (0.6,2.6) node[above] {\footnotesize $a$} circle  (2.2pt);
\filldraw (3.6,2.6) node[above] {\footnotesize $b$} circle  (2.2pt);
\filldraw (1.1,1.6) node[left] {\footnotesize $c$} circle  (2.2pt);
\filldraw (4.1,1.6) node[right] {\footnotesize $d$} circle  (2.2pt);
\filldraw (0.6,0.6) node[left] {\footnotesize $e$} circle  (2.2pt);
\filldraw (3.6,0.6) node[right] {\footnotesize $f$} circle  (2.2pt);
\filldraw (0,2) node[left] {\footnotesize $p$} circle  (2.2pt);
\filldraw (3,2) node[right] {\footnotesize $q$} circle  (2.2pt);
\filldraw (0.5,1) node[left] {\footnotesize $r$} circle  (2.2pt);
\filldraw (3.5,1) node[right] {\footnotesize $s$} circle  (2.2pt);
\filldraw (0,0) node[left] {\footnotesize $t$} circle  (2.2pt);
\filldraw (3,0) node[right] {\footnotesize $u$} circle  (2.2pt);

  \node at (1.75,-0.75) {$N_2$};  
\end{tikzpicture}
\caption{The rank-$4$ paving matroids in Example
  \ref{ex:4}.}\label{fig:ex4}
\end{figure}
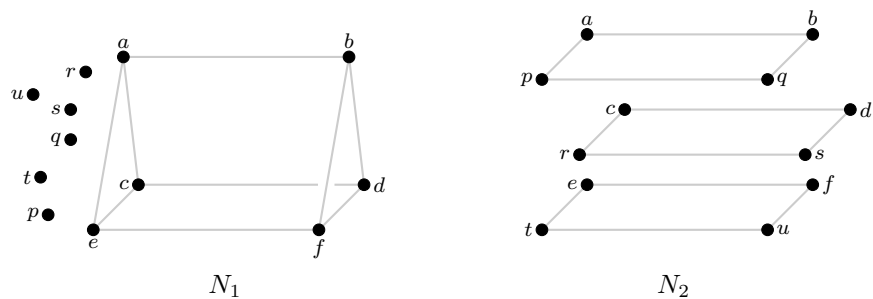

\begin{example}\label{ex:4}
  Let $N_1$ and $N_2$ be the rank-$4$ paving matroids on
  $\{a,b,c,d,e,f,q,p,r,s,t,u\}$ where
  \begin{itemize}
  \item the dependent hyperplanes of $N_1$ are $\{a,b,c,d\}$,
    $\{a,b,e,f\}$, and $\{c,d,e,f\}$, and
  \item the dependent hyperplanes of $N_2$ are $\{a,b,p,q\}$,
    $\{c,d,r,s\}$, and $\{e,f,t,u\}$.
  \end{itemize}
  These paving matroids are shown in Figure \ref{fig:ex4}. 
  We partition $\mathcal{H}_1-\mathcal{H}_2$ into
  \begin{itemize}
  \item
    $A_1= \bigl\{\{a,b,e,f\},\{a,b,p\},\{a,b,q\},\{a,p,q\},\{b,p,q\}
    \bigr\}$,
  \item
    $A_2= \bigl\{\{a,b,c,d\},\{c,d,r\},\{c,d,s\},\{c,r,s\},\{d,r,s\}
    \bigr\}$,
  \item
    $A_3= \bigl\{\{c,d,e,f\},\{e,f,t\},\{e,f,u\},\{e,t,u\},\{f,t,u\}
    \bigr\}$,
  \end{itemize}
  and $\mathcal{H}_1-\mathcal{H}_2$ into
  \begin{itemize}
  \item
    $B_1=
    \bigl\{\{a,b,p,q\},\{a,b,e\},\{a,b,f\},\{a,e,f\},\{b,e,f\}\bigr\}$,
  \item
    $B_2= \bigl\{\{c,d,r,s\},\{c,d,a\},\{c,d,b\},\{c,a,b\},\{d,a,b\}
    \bigr\}$,
  \item
    $B_3= \bigl\{\{e,f,t,u\},\{e,f,c\},\{e,f,d\},\{e,c,d\},\{f,c,d\}
    \bigr\}$.
  \end{itemize}
  Three bijections that have the properties in Theorem
  \ref{thm:generalpaving}, given as products of cycles, are
  $\alpha_1=(e,p)(f,q)$, $\alpha_2=(a,r)(b,s)$, and
  $\alpha_3=(c,t)(d,u)$.

  The lattices of flats, $L_1$ and $L_2$, of $N_1$ and $N_2$ have
  $291$ elements, and the isomorphisms
  $\phi_1:L_1\to \mathcal{Z}(M_1)$ and
  $\phi_2:L_2\to \mathcal{Z}(M_2)$ must, for instance, assign to
  $\{e\}$ and $\{p\}$ cyclic flats of the same size and rank, and
  likewise for $\{a,e\}$ and $\{a,p\}$.
\end{example}  

\begin{center}
 \textsc{Acknowledgments}
\end{center}

\vspace{3pt}

The author thanks both referees for their careful reading of the
manuscript and their thoughtful comments.


\begin{thebibliography}{99}

\bibitem{asymp} N.~Bansal, R.A.  Pendavingh, and J.~van der Pol, On
  the number of matroids, \emph{Combinatorica} \textbf{35} (2015)
  253--277.
  
\bibitem{catdata} J.~Bonin and J.P.S.~Kung, The
  $\mathcal{G}$-invariant and catenary data of a matroid, \emph{Adv.\
    in Appl.\ Math.}\ \textbf{94} (2018) 39--70.

\bibitem{cone} J.~Bonin and K.~Long, The free $m$-cone of a matroid
  and its $\mathcal{G}$-invariant, submitted.
	
\bibitem{cyclic} J.\ Bonin and A.\ de Mier, The lattice of cyclic
  flats of a matroid, \emph{Ann.\ Comb.}\ \textbf{12} (2008) 155--170.
  
\bibitem{survey} J.~Bonin and A.~de~Mier, Tutte uniqueness and Tutte
  equivalence, In \cite{handbook} (in press; 2022).

\bibitem{Tutte} T.~Brylawski and J.~Oxley, The Tutte Polynomial and
  Its Applications, In N. White (ed), \emph{Matroid Applications}
  (Encyclopedia of Mathematics and its Applications,
  pp. 123--225). Cambridge: Cambridge University Press (1992).
	
\bibitem{G-inv} H.~Derksen, Symmetric and quasi-symmetric functions
  associated to polymatroids, \emph{J.\ Algebraic Combin.}\
  \textbf{30} (2009) 43--86.
	
\bibitem{valuative} H.~Derksen and A.~Fink, Valuative invariants for
  polymatroids, \emph{Adv.\ Math.}\ \textbf{225} (2010) 1840--1892
	
\bibitem{config} J.\ Eberhardt, Computing the Tutte polynomial of a
  matroid from its lattice of cyclic flats, \emph{Electron.\ J.\
    Combin.}\ \textbf{21} (2014) Paper 3.47, 12 pp.

\bibitem{handbook} J.A.~Ellis-Monaghan and I.~Moffatt (eds),
  \emph{Handbook of the Tutte Polynomial and Related Topics}.  (in
  press) CRC Press/Taylor \& Francis (2022).

\bibitem{oxley} J.~Oxley, \emph{Matroid Theory}, second edition
  (Oxford University Press, Oxford, 2011).

\bibitem{sims} J.A.\ Sims, \emph{Some Problems in Matroid Theory},
  (Ph.D. Dissertation, Linacre College, Oxford University, Oxford,
  1980).
  
\end{thebibliography}
\end{document}